\documentclass[twoside,12pt]{article}
\usepackage[usenames,dvipsnames]{color}

\usepackage{amsmath}
\usepackage{amssymb}
\usepackage{mathrsfs}
\usepackage{verbatim}
\usepackage{cite}

\definecolor{rosso}{rgb}{0.8,0,0}

\def\Sign{\mathop{\rm Sign}}

\def\interna{\mathop{\rm int}}

\def\Signe{\mathop{\rm Sign_{\varepsilon}}}

\title{On a class of conserved phase field systems \\
with a maximal monotone perturbation}

\author{Michele Colturato\\
Dipartimento di Matematica, Universit\`a degli Studi di Pavia\\
Via Ferrata~1, 27100 Pavia, Italy\\
E-mail: \texttt{michele.colturato01@universitadipavia.it}}
 
\date{}

\newcommand\testopari{\sc Michele Colturato}
\newcommand\testodispari{\sc On a class of conserved phase field systems }
\newtheorem{Teo1-Esistenza}[def1]{Theorem}
\newtheorem{Teo11-Regolare}[def1]{Theorem}
\newtheorem{Teo2-Unicita-Dipendenzacontinua}[def1]{Theorem}
\newtheorem{Teo3-Sliding}[def1]{Theorem}
\newtheorem{Lemma-Sliding}[def1]{Lemma}

\begin{document}
\maketitle

\begin{abstract}
We prove existence and regularity for the solutions to a Cahn--Hilliard system describing the phenomenon of phase separation
for a material contained in a bounded and regular domain. Since the first equation of the system
is perturbed by the presence of an additional maximal monotone operator,
we show our results using suitable regularization
of the nonlinearities of the problem and performing some a priori estimates which allow
us to pass to the limit thanks to compactness and monotonicity arguments.
Next, under further assumptions, we deduce a continuous dependence estimate 
whence the uniqueness property is also achieved.
Then, we consider the related sliding mode control (SMC) problem   
and show that the chosen SMC law forces a suitable linear combination of the temperature and the phase to
reach a given (space-dependent) value within finite time.

\vspace{2mm}

\vspace{2mm}
\noindent \textbf{AMS (MOS) subject clas\-si\-fi\-ca\-tion:} 35K61, 35K25, 35D30, 34H05, 80A22.

\end{abstract}

%%%%% Section 1. %%%%% 
\section{Introduction}
The  Cahn--Hilliard e quation, originally
introduce d in \cite{CH58} and first studie d mathe matically in the  se minal pape r 
\cite{EZ86}, yie lds a de scription of the  e volution phe nome non of
the  solid--solid phase  se paration. 
In ge ne ral, an e volution proce ss goe s on diffusive ly.
Howe ve r, the  phe nome non of the  solid--solid phase  se paration doe s not se e m to follow
this structure : e .g., whe n a binary alloy is coole d down sufficie ntly, e ach phase  conce ntrate s and
the  mate rial quickly be come s inhomoge ne ous, forming a fine -graine d structure  in which
e ach of the  two compone nts appe ars more  or le ss alte rnative ly (se e , e.g., \cite{KN96}). 
The Cahn--Hilliard e quation is a ce le brate d mode l which
de scribe s this proce ss (usually known as spinodal de composition)
by the  simple  frame work of partial diffe re ntial e quations.
The  mathe matical lite rature  conce rning this proble m is rathe r vast.
Le t us quote  \cite{eq1, CF15bbbbbb, GMS09, preamboloiniziale, Kub12, eq2, eq3, eq4, eq5} and also re fe r to \cite{CF15a} in which  a
force d mass constraint on the  boundary is conside re d. 

In the  pre se nt contribution, we  conside r the  following Cahn--Hilliard syste m 
pe rturbe d by the  pre se nce  of an additional maximal monotone  nonline arity: 
\begin{equation} \label{iniz-1}
\partial_t(\chi + \ell \varphi) - \Delta \chi + \zeta = f \ \ \textrm{a.e. in $Q:=  \Omega \times (0,T) $,} 
\end{equation}
\begin{equation} \label{iniz-2}
\partial_t \varphi - \Delta \mu= 0  \ \ \textrm{a.e. in $Q$,}
\end{equation}
\begin{equation} \label{iniz-3}
\mu= - \nu \Delta \varphi + \xi + \pi(\varphi) - \gamma \chi \ \ \textrm{a.e. in $Q$,} 
\end{equation}
\begin{equation} \label{iniz-4}
\zeta(t) \in A(a \chi(t) + b \varphi(t) - \kappa^*) \ \ \textrm{for a.e. $t \in (0,T)$,}  
\end{equation}
\begin{equation} \label{iniz-5}
\xi \in \beta(\varphi) \ \ \textrm{a.e. in $Q$,}
\end{equation}
whe re  $\Omega \subseteq \mathbb{R}^3$ is an open, bounde d, conne cte d subse t of class $C^1$,
$T$~is some  final time ,
$\chi$, $\varphi$ and $\mu$ de  note   the   te  mpe  rature  , the  orde r parame te r and the  che mical pote ntial, re spe ctive ly.
We  point out that he re  $\chi$ doe s not re pre se nt the  absolute  te mpe rature ,
but it is re late d to it by 
\begin{equation} \label{temperaturaassoluta}
\chi = \Theta - \Theta_c,
\end{equation}
whe re $\Theta_c$ denote s a critical te mpe rature .
More ove r, $\kappa^*$ is a function
in $H^2(\Omega)$ with null outward normal de rivative  on the  boundary of $\Omega$,
$f$~is a source term and $a$, $b$, $\ell$, $\gamma$ are constants.
In particular, let $\ell$ and $\gamma$ be positive.
The above syste m is comple me nte d by
homoge ne ous Ne umann boundary conditions for both $\chi$ and~$\varphi$, that is,
\begin{equation} \label{iniz-6}
\partial_{\bf{n}} \chi = \partial_{\bf{n}} \varphi = \partial_{\bf{n}} \mu = 0 \quad \textrm{on $\Sigma:= \Gamma \times (0,T)$},
\end{equation}
where $\Gamma$ is the boundary of $\Omega$ and $\partial_{\bf{n}}$ is the outward normal derivative,
and by the initial conditions 
\begin{equation} \label{iniz-7}
\chi(0) = \chi_0, \quad \quad \varphi(0) = \varphi_0 \quad \textrm{in $\Omega$}.
\end{equation}
The term $\xi + \pi(\varphi)$, appearing in \eqref{iniz-3}, represents the derivative of 
the pote ntial $\mathcal{W}$ associate d with the  phase  configuration. In the 
lite rature  (se e, e.g., \cite{CoGiMaRo, EllZheng}), $\mathcal{W}$ is frequently assume d to be  a double -we ll pote ntial. 
More  ge ne rally, $\mathcal{W}$ can be  de fine d as the  sum $\mathcal{W}= \widehat{\beta} + \widehat{\pi}$,
where $\widehat{\beta}: \mathbb{R} \rightarrow [0, + \infty]$ is a proper, l.s.c. and conve x function
and $ \widehat{\pi}:  \mathbb{R} \rightarrow \mathbb{R}$ is a function in $C^1(\mathbb{R})$ such that $\pi:= \widehat{\pi}'$ is Lipschitz continuous. 
Due  to the  prope rtie s of $\widehat{\beta}$, the  subdiffe re ntial $\partial \widehat{\beta} =: \beta$ is we ll de fine d and is a maximal monotone  graph.
In our proble m a maximal monotone  ope rator $A: \ H:=L^2(\Omega) \rightarrow  2^H$ also appe ars. We  assume  that
$0 \in A(0)$ and $ \| v \|_H \leq C (1 + \| x \|_H )$  for all $x \in H$, $v \in A x$,
for some constant $C>0$.  For a compre he nsive  pre se ntation of the  the ory of 
maximal monotone  ope rators, we  re fer, e.g., to \cite{Bar10, Bre73, Show}.

Le t us spe nd some  words about the  the rmodynamic de rivation of the system
\eqref{iniz-1}--\eqref{iniz-5}. We move from the following expression for the local free energy:
\begin{equation}
\Psi (\Theta , \varphi)= c_0 \Theta (1- \ln \Theta) - \gamma_0 (\Theta - \Theta_c) \varphi
+ \Theta \bigg( - \mu \varphi + \mathcal{W}(\varphi) + \frac{\nu}{2}|\nabla \varphi|^2 \bigg),
\end{equation}
where $\Theta$ denotes the absolute temperature, as in \eqref{temperaturaassoluta}.
The quantity 
\begin{equation}
e = \chi + \ell \varphi
\end{equation}
appearing (under time derivative) in \eqref{iniz-1} represents a rescaled internal energy 
of the system, since a standard thermodynamic relation postulates that
\begin{equation}
\mathcal{U} = \Psi + \Theta \mathcal{S},
\end{equation}
where $\mathcal{U}$ is the internal energy and $\mathcal{S}:= - \frac{\partial \Psi}{\partial \Theta}$ denotes the entropy.
Then, by an easy computation, we find out that
\begin{align} \nonumber
&\mathcal{U} = \ 
c_0 \Theta (1- \ln \Theta) - \gamma_0 (\Theta - \Theta_c) \varphi
+ \Theta \bigg( - \mu \varphi + \mathcal{W}(\varphi) + \frac{\nu}{2}|\nabla \varphi|^2 \bigg) \nonumber\\
&{}- \Theta \bigg( c_0(1- \ln \Theta) -c_0 - \gamma_0 \varphi - \mu \varphi + \mathcal{W}(\varphi) + \frac{\nu}{2}|\nabla \varphi|^2 \bigg)
= c_0 \Theta + \gamma_0 \Theta_c \varphi,
\end{align}
so that, by adding the constant $-c_0 \Theta_c$
and dividing by $c_0$, from \eqref{temperaturaassoluta} we exactly have that
\begin{equation}
e = \chi + \ell \varphi, \quad \quad \textrm{with} \quad \quad \ell = \frac{\gamma_0 \Theta_c}{c_0}.
\end{equation}
On the other hand, by considering the total free energy
\begin{equation}
\mathcal{F} (\Theta , \varphi)= \int_{\Omega} \Psi (\Theta , \varphi)
\end{equation}
and taking the  variational de rivative  $\frac{\delta \mathcal{F} }{\delta \varphi}$,
we actually recover the phase equation \eqref{iniz-3}. Indeed, by a standard computation, we infer that 
\begin{equation}
- \gamma_0 (\Theta - \Theta_c) \varphi
+ \Theta \big( - \mu + \mathcal{W}'(\varphi) - \nu \Delta \varphi \big)
- \nu \nabla \Theta \cdot \nabla \varphi = 0.
\end{equation} 
Dividing by $\Theta$, observing that (cf. \eqref{temperaturaassoluta})
\begin{equation}
- \gamma_0 \bigg(1 - \frac{\Theta_c}{\Theta} \bigg) \cong  \frac{\gamma_0}{\Theta_c} \chi,
\end{equation} 
thanks to a first orde r line arization, and ne gle cting the  highe r orde r term $-\frac{\nu}{\Theta} \nabla \Theta \cdot \nabla \varphi$,
we finally obtain \eqref{iniz-3} with $\gamma \cong \frac{\gamma_0}{\Theta_c}$.

As usual for Cahn--Hilliard syste m, in the  Proble m $(P)$ state d by \eqref{iniz-1}--\eqref{iniz-7} the  inte gral me an value  of $\varphi(t)$ 
re mains constant during the  whole  e volution. Inde ed, fixing 
an arbitrary $t \in (0,T)$ and inte grating \eqref{iniz-2} over $\Omega$, 
we infe r that
\begin{equation} \label{utileinregolaritaitti-intro}
\frac{d}{dt} \int_{\Omega} \varphi(t) = 0,
\end{equation}
whe nce  it imme diate ly follows that
\begin{equation} \label{integraleconservato11-intro}
m(\varphi(t)) := \frac{1}{|\Omega|} \int_{\Omega} \varphi(t) = \frac{1}{|\Omega|} \int_{\Omega} \varphi_0 
\quad \quad \textrm{for every $t \in (0,T)$}.
\end{equation}
We  also obse rve  that syste m $(P)$  is a fourth-orde r proble m
constructe d as the  conse rve d ve rsion of the  following phase -fie ld syste m:
\begin{equation} \label{originale1111}
\partial_t(\chi + \ell \varphi) -k \Delta \chi + \zeta = f \ \ \textrm{a.e. in $Q$,} 
\end{equation}
\begin{equation} \label{originale21111}
\partial_t \varphi - \nu \Delta \varphi + \xi + \pi(\varphi) = \gamma \chi  \ \ \textrm{a.e. in $Q$,}
\end{equation}
\begin{equation} 
\zeta(t) \in A(\chi(t) + \alpha \varphi(t) - \kappa^*) \ \textrm{for a.e. $t \in (0,T)$,}  
\end{equation}
\begin{equation} 
\xi \in \beta(\varphi) \ \textrm{a.e. in $Q$,}
\end{equation}
\begin{equation}
\partial_{\bf{n}}\chi = 0, \quad \quad \partial_{\bf{n}}\varphi = 0 \quad \textrm{on $\Sigma$},
\end{equation}
\begin{equation} \label{originale51111}
\chi(0)=\chi_0, \quad \quad \varphi(0)=\varphi_0  \quad  \textrm{in $\Omega$,}
\end{equation} 
whe re $k$ and $\alpha$ are  positive  coe fficie nts.
Phase -fie ld syste ms have  be e n wide ly studie d in the  lite rature . 
We  re fe r, without any sake  of comple te ne ss, e .g., to \cite{BrokSpr, Cag, CLS99, EllZheng, GraPetSch, KenmNiez, KN96, Lau, sprekzen} 
and re fe re nce s the re in for the  we ll-pose dne ss and long-time  be havior re sults.
In particular, the  above  syste m \eqref{originale1111}--\eqref{originale51111}
has be e n thoroughly discusse d in \cite {Michele},
whe re  e xiste nce  and re gularity of the  solutions is prove d
and, unde r furthe r assumptions, unique ne ss and continuous de pe nde nce  on the  initial data are  de duce d.

In the  pape r, we  first show the  e xiste nce  of solutions for Proble m $(P)$ (se e \eqref{iniz-1}--\eqref{iniz-7}).
In order to carry out this purpose, we  conside r the  approximating proble m ($P_{\varepsilon}$),
obtaine  d from ($P$) by approximating $A$ and $\beta$ by the ir Yosida
re gularizations. In pe rforming our uniform e stimate s
we  ofte n re fe r to \cite{CF15200009}, whe re   the  authors propose  the  study
of a nonline ar diffusion proble m as an asymptotic limit of a particular Cahn--Hilliard syste m.
The n, we  pass to the  limit as $\varepsilon \searrow 0$
and show that some  limit of a subse que nce  of solutions for ($P_{\varepsilon}$) yie lds a solution of ($P$).
Ne xt, we le t $a\ell = b$ which is, in some  se nse , a physical re striction 
since  the  argume nt of the  variable  in the  ope rator $A$ is thus proportional to the  inte rnal e ne rgy of the  syste m.
We  also write  Proble m $(P)$ for two diffe re nt se ts of data $f_i$, $\kappa^*_i$, $\chi_{0_{i}}$ and $\varphi_{0_{i}}$, $i=1,2$.
By suitably pe rforming contracting e stimate s for the  diffe re nce  of the  corre sponding solutions,
we  de duce  the  continuous de pe nde nce  re sult whe nce  the  unique ne ss prope rty is also achie ve d.

The  se cond part of the  pape r is de vote d to the  sliding mode  control (SMC) proble m. 
He nce , the  main ide a be hind this sche me  is first to ide ntify a manifold of lowe r dime nsion (calle d the  sliding manifold) 
whe re  the  control goal is fulfille d and such that the  original syste m re stricte d to this sliding manifold has a
de sire d be havior, and the n to act on the  syste m through the  control in orde r to constrain
the  e volution on it, that is, to de sign a SMC-law that force s the  traje ctorie s of the  syste m
to re ach the  sliding surface  and maintains the m on it (se e , e .g., \cite{I76, O00}).
The  main advantage  of sliding mode  control is that it allows the  se paration of the 
motion of the  ove rall syste m in inde pe nde nt partial compone nts of lowe r dime nsions, and
conse que ntly it re duce s the  comple xity of the  control proble m. 
In particular, we  prove  the  e xiste nce  of sliding mode s for the  solutions of our syste m $(P)$
for a suitable  choice  of the  ope rator $A$ and the  coe fficie nts $a$ and $b$.
We  take  $a= 1$, $b= \ell$ and $A = \rho \Sign$, whe re 
$\rho$ is a positive  coe fficie nt and $\Sign : H \longrightarrow 2^H$ is a maximal monotone  ope rator de fine d as
$\Sign(v) = \frac{v}{\| v \|}$, if $v \neq 0$ and $\Sign(0) = B_1(0)$, if $v = 0$ (here,
$B_1(0)$ denotes the close d unit ball of $H$).
Thus we  pre scribe  a state -fe e dback control law acting on  the  re scale d inte rnal e ne rgy $(\chi + \ell \varphi)$
of the   syste  m in orde  r that the   dynamics of the   syste  m modifie  d in this
way force  s the   value  $(\chi(t) + \ell \varphi(t))$ to reach a manifold of the  phase  space  
in a finite  time  and the n lie  the re  with a sliding mode  (cf. \cite{BaCoGiMaRo, collisprek}). 

Conce rning the  study of optimal control proble ms for phase -fie ld syste ms,
we  quote \cite{CGM, CoGiMaRo, CoGiMaRorrrrrrrrr, HKKY}.
Re ce nt inve stigations have  be e n also addre sse d to the  optimal control proble m 
for Cahn-Hilliard syste ms: le t us me ntion \cite{SC1, SC2, SC3, SC4, control1, eq2}.
We  also re fe r to \cite{control4, control5} which
deals with the  conve ctive  Cahn--Hilliard e quation, 
and to \cite{control7, control8}, whe re
some  discre tize d ve rsions of the  ge ne ral Cahn--Hilliard syste ms are  studie d.

In the  pre se nt contribution, assuming  $a= 1$, $b=\ell$ and $A = \rho \Sign$
in \eqref{iniz-1}--\eqref{iniz-7}, we  prove  the  e xiste nce  of sliding mode s for Proble m $(P)$
by ide ntifying $\rho^* > 0$ such that the  following prope rty is fulfille d: 
for e ve ry $\rho >\rho^*$, the re  e xists a solution $(\chi,\varphi,\mu)$
to Proble m $(P)$ and a time  $T^*$ such that, for e ve ry $t \in [T^*,T]$
\begin{equation}
\chi(t) + \ell \varphi(t) = \kappa^* \quad \quad \textrm{a.e. in $\Omega$.}
\end{equation}
It is curious and inte re sting that we  are  able  to handle 
a fe e dback law and prove  the  me ntione d prope rty just for the  inte rnal e ne rgy
of the  syste m, which is a spe cial line ar combination of the  variable s $\chi$
and $\varphi$. Howe ve r, for a discussion of the  SMC laws, line ar and nonline ar,
that can be  conside re d for phase  fie ld syste ms, we  re fe r to the  Introduction of \cite{BaCoGiMaRo}.

The  pape r is organize d  as follows.
In Se ction 1, we  list our assumptions, state  the  proble m in a pre cise  form
and pre se nt our re sults.
The  ne xt se ctions are  de vote d to the  corre sponding proofs:
Se ction 3--6 de al with e xiste nce  and re gularity,
while  unique ne ss and continuous de pe nde nce  are  prove d in Se ction 7.
In Se ction 8, we  show the  e xiste nce  of sliding mode s.

%%%%% Se ction 2. %%%%%
\section{Main results}
\setcounter{equation}{0}
In this section, we state the main results. 

%%%%% Section 2.1. %%%%%
\subsection{Preliminary assumptions}
We  assume  $\Omega \subseteq \mathbb{R}^3$ to be ope n, bounde d, conne cte d, of class $C^1$ and we write $| \Omega |$ for its Lebesgue measure.
Moreover, $\Gamma$ and $\partial_{\bf{n}}$ still stand for the boundary of $\Omega$ and the outward normal derivative, respectively.
Given a finite final time $T > 0$, for every $t \in (0,T]$ we set
\begin{equation}
Q_t = (0,t) \times \Omega, \ \ Q = Q_T,  
\end{equation}
\begin{equation}
\Sigma_t = (0,t) \times \Gamma, \ \ \Sigma = \Sigma_T. 
\end{equation}
In the following, we set for brevity:
\begin{equation}\label{W} 
H = L^2(\Omega), \ \ \ \ V = H^1(\Omega), \ \ \ \ V_0 = H^1_0(\Omega), \ \ \ \
W = \{ u \in H^2(\Omega): \ \partial_{\bf{n}} u = 0 \ \textrm{on} \ \partial \Omega \},
\end{equation}
with usual norms $\| \cdot \|_{H}$, $\|\cdot \|_{V}$
and inner products $(\cdot,\cdot )_{H}$, $(\cdot ,\cdot )_{V}$, respectively. 
The symbol $V^*$ denotes the dual space of $V$ while the pair $\langle \cdot ,\cdot \rangle _{V^*, V}$ represents
the duality pairing between $V^*$ and $V$. Moreover, we identify $H$ with its dual space.

\subsection{Operators}
In this subsection we describe the operators appearing in the problem under study.
\paragraph{The operator $m$.} 
We consider the operator $m:V^* \rightarrow \mathbb{R}$ defined by
\begin{equation} \label{m1}
	m(z^*):=\frac{1}{|\Omega |}
	\langle z^*, 1 \rangle _{V^*,V}
	\quad \hbox{for~all~} z^* \in V^*.
\end{equation}
We observe that, if $z^* \in H$, then
\begin{equation} \label{m2}
	m(z^*)=\frac{1}{|\Omega |}
	\int_{\Omega }^{}z^* dx.
\end{equation}
\paragraph{The double-well potential $\mathcal{W}$.}
We introduce the double-well potential $\mathcal{W}$ as the sum 
\begin{equation} 
\mathcal{W}= \widehat{\beta} + \widehat{\pi},
\end{equation}
where
	\begin{equation} \label{beta}
\widehat{\beta}: \mathbb{R} \longrightarrow [0, + \infty] \textrm{ is proper, l.s.c. and convex with $\widehat{\beta}(0)=0$,} 
		\end{equation}
		\begin{equation} \label{pi}
\widehat{\pi}:  \mathbb{R} \rightarrow \mathbb{R}, \ \textrm{$\widehat{\pi} \in C^1(\mathbb{R})$ with $ \pi := \widehat{\pi}'$ Lipschitz continuous.}
	\end{equation}
Since $\widehat{\beta}$ is proper, l.s.c. and convex, the subdifferential $\beta:= \partial \widehat{\beta} $ is well defined. We denote by 
$D(\beta)$ and $D(\widehat{\beta})$ the effective domains of $\beta$ and $\widehat{\beta}$, respectively, and also assume that $\interna(D(\beta)) \neq \emptyset$. Thanks to these assumptions, $\beta$ is a maximal monotone graph.	Moreover, as $\widehat{\beta}$ 
takes its minimum in $0$, we have that $0 \in \beta(0)$. 

\paragraph{The operator $\mathcal{B}$.}
We introduce the operator $\mathcal{B}$ induced by $\beta$ on $L^2(Q)$ in the following way:
\begin{equation} \label{betagrande1}
\mathcal{B}: L^2(Q) \longrightarrow L^2(Q)
\end{equation}
\begin{equation} \label{betagrande2}
\xi  \in \mathcal{B}(\varphi) \Longleftrightarrow \xi(x,t) \in \beta(\varphi(x,t)) \quad \textrm{ for a.e. $(x,t) \in Q$.}
\end{equation}
We notice that
\begin{equation} 
\beta= \partial \widehat{\beta},  \quad \quad \quad \quad \mathcal{B} = \partial \Phi, 
\end{equation}
where
\begin{equation} 
\Phi: \  L^2(Q) \longrightarrow (-\infty , + \infty]
\end{equation}
\begin{equation}
\Phi(u)=
\left\{ \begin{array}{ll}
\int_Q{\widehat{\beta}(u)} & \textrm{if $u \in L^2(Q)$ and $\widehat{\beta}(u) \in L^1(Q)$}, \\
+ \infty & \textrm{elsewhere, with  $u \in L^2(Q)$.}
\end{array}
\right.
\end{equation}

\paragraph{The operator $A$.}
We consider the maximal monotone operator
\begin{equation} \label{A1}
A: H \longrightarrow  H.
\end{equation}
We assume that
\begin{equation} \label{A2}
0 \in A(0)
\end{equation}
and that there exists a constant $C_A>0$ such that
\begin{equation} \label{stimaA}
\| v \|_H \leq C_A (1 + \| \kappa \|_H ) \quad \textrm{for every $\kappa \in H$, $v \in A\kappa$}.
\end{equation}

\paragraph{The operator $\mathcal{A}$.}
We introduce the operator $\mathcal{A}$ induced by $A$ on $L^2(0,T;H)$ in the following way
\begin{equation} \label{Agrande1}
\mathcal{A}: L^2(0,T;H) \longrightarrow L^2(0,T;H)
\end{equation}
\begin{equation} \label{Agrande2}
\zeta  \in \mathcal{A}(\kappa) \Longleftrightarrow \zeta(t) \in A(\kappa(t)) \quad \textrm{ for a.e. $t \in (0,T)$.}
\end{equation}
We notice that also $\mathcal{A}$ is a maximal monotone operator.

\paragraph{The operator $\Sign$.}
An example of maximal monotone operator $A$ which satisfies \eqref{A1}--\eqref{stimaA}
is the operator 
\begin{equation} \label{Sign}
\Sign: \quad \quad H \longrightarrow 2^H
\end{equation}
\begin{equation} \label{Sign3333333}
\Sign(v) = 
\left\{ \begin{array}{ll}
\frac{v}{\| v \|} 								& \textrm{if $v \neq 0$}, \\
B_1(0)													& \textrm{if $ v =0  $,}
\end{array}
\right.
\end{equation}
where $B_1(0)$ is the closed unit ball of $H$.
$\Sign$ is the subdifferential of the map $ \| \cdot \| : H \rightarrow  \mathbb{R} $ and
is a maximal monotone operator on $H$ which satisfies \eqref{A2}--\eqref{stimaA}.

\paragraph{The operator $\mathcal{N}$.}
We also consider the operator 
\begin{equation} \label{operatoreN1}
\mathcal{N}:D({\mathcal N}) \subseteq V^* \to V, \quad
\end{equation}
defined on its domain 
\begin{equation} \label{operatoreNNNN1}
D({\mathcal N}):=\{w \in V^* : m(w^*)=0 \}.
\end{equation}
For every $w^* \in D({\mathcal N})$,  we define
$w={\mathcal N}w^*$ if $w \in V$, $m(w)=0$ and $w$ is a solution 
of the following variational equation
\begin{equation} \label{operatoreN2}
	\int_{\Omega }^{}\nabla w\cdot \nabla z dx = \langle w^*,z \rangle _{V^*,V}
	\quad \hbox{for~all~} z \in V.
\end{equation}
If $w^* \in D({\mathcal N}) \cap H$, then $w$ is the unique solution to the elliptic problem
\begin{equation} \label{operatoreN3}
	\begin{cases}
	\displaystyle - \Delta w = w^* & \hbox{a.e.\ in~} \Omega, \vspace{1mm}\\
	\partial _{\nu }  w = 0 & \hbox{a.e.\ in~} \Gamma, \vspace{1mm}\\
	m(w)  = 0 . &
	\end{cases} 
\end{equation}
We observe that, due to elliptic regularity, $w\in W$. Moreover, for every
$v^*,w^* \in D({\mathcal N})$, $v={\mathcal N}v^* $ and $w={\mathcal N}w^* $ we have that
\begin{align}
	\langle w^*, {\mathcal N}v^* \rangle _{V^*,V} 
	& = \langle w^*, v \rangle _{V^*,V} \nonumber = \int_{\Omega } \nabla w \cdot \nabla v dx \nonumber \\
	& = \langle v^*, w \rangle _{V^*,V} \nonumber = \langle v^*, {\mathcal N}w^* \rangle _{V^*,V}.
\end{align}
Consequently, by defining
\begin{equation} 
	\|w^*\|_{V^*}^2:=\bigl \|\nabla {\mathcal N} \bigl( w^*-m(w^*) \bigr )\bigr \|_{H^3}^2
	+ \bigl| m(w^*) \bigr|^2 
	\quad \hbox{for~all~}w^* \in V^*,
	\label{norm}
\end{equation}
it turns out that
$\|\cdot \|_{V^*}$ is a norm in $V^*$. 

\subsection{Setting of the problem and results}
Now, we describe the state system. We assume
\begin{equation} \label{parametri}
\ell, \ \nu, \ \gamma \in (0, + \infty), \quad \quad a, \ b \in \mathbb{R},
\end{equation}
\begin{equation} \label{f}
f \in L^2(0,T, H),
\end{equation}
\begin{equation} \label{0}
\kappa^{*} \in W, \quad \chi_0 \in H, \quad \varphi_0 \in V, 
\quad  \widehat{\beta}(\varphi_0) \in L^1(\Omega), \quad  m(\varphi_0) =: m_0 \in \interna(D(\beta)).
\end{equation}
We look for a triplet $(\chi,\varphi,\mu)$ 
satisfying at least the regularity requirements
\begin{equation} \label{regolaritaesistenza1}
\chi  \in H^1(0,T;V^*) \cap L^\infty (0,T;H) \cap  L^2(0,T;V), 
\end{equation}
\begin{equation} \label{regolaritaesistenza2}
\varphi    \in H^1(0,T;V^*) \cap L^\infty (0,T;V) \cap L^2(0,T;W), 
\end{equation}
\begin{equation} \label{regolaritaesistenza3}
\mu        \in L^2(0,T;V),
\end{equation}
and solving the Problem ($P$), that is,
\begin{equation} \label{iniziale1}
\partial_t(\chi + \ell \varphi) - \Delta \chi + \zeta = f \ \ \textrm{a.e. in $Q$,} 
\end{equation}
\begin{equation} \label{iniziale2}
\partial_t \varphi - \Delta \mu= 0  \ \ \textrm{a.e. in $Q$,}
\end{equation}
\begin{equation} \label{iniziale3}
\mu= - \nu \Delta \varphi + \xi + \pi(\varphi) - \gamma \chi \ \ \textrm{a.e. in $Q$,} 
\end{equation}
\begin{equation} \label{iniziale4}
\zeta(t) \in A(a \chi(t) + b \varphi(t) - \kappa^*) \ \textrm{for a.e. $t \in (0,T)$,}  
\end{equation}
\begin{equation} \label{iniziale5}
\xi \in \beta(\varphi) \ \textrm{a.e. in $Q$,}
\end{equation}
\begin{equation} \label{iniziale6}
\partial_{\bf{n}} \chi = \partial_{\bf{n}} \varphi = \partial_{\bf{n}} \mu = 0 \quad \textrm{on $\Sigma$,}
\end{equation}
\begin{equation} \label{iniziale7}
\chi(0) = \chi_0, \quad \quad \varphi(0) = \varphi_0 \quad \textrm{in $\Omega$}.
\end{equation}

\begin{Teo1-Esistenza}[Existence] \label{Teorema-esistenza}
Assume \eqref{beta}--\eqref{pi}, \eqref{A1}--\eqref{stimaA} and \eqref{parametri}--\eqref{0}. 
Then Problem $(P)$ (see \eqref{iniziale1}--\eqref{iniziale7})
has at least one solution $(\chi,\varphi,\mu)$ satisfying \eqref{regolaritaesistenza1}--\eqref{regolaritaesistenza3}. 
\end{Teo1-Esistenza}

\begin{Teo11-Regolare}[Regularity] \label{Regolar}
Assume \eqref{beta}--\eqref{pi}, \eqref{A1}--\eqref{stimaA}, \eqref{parametri}--\eqref{f}, 
\begin{equation} \label{0reg}
\kappa^{*} \in W, \quad  \chi_0 \in V, \quad  \varphi_0 \in W,
\quad \beta^0(\varphi_0) \in H, \quad m_0 \in \interna(D(\beta))
\end{equation}
and that there exists $\varepsilon_0 \in (0,1]$
such that
\begin{equation} \label{tecnicadistima1}
\| -\nu \Delta \varphi_0 + \beta_{\varepsilon}(\varphi_0) + \pi(\varphi_0) - \gamma \chi_0 \|_V \leq c \quad \quad \textrm{for every $\varepsilon \in (0, \varepsilon_0]$},
\end{equation}
for some positive constant $c$, then the solution $(\chi,\varphi,\mu)$ satisfies 
\begin{equation} \label{regolaritaesistenza1reg}
\chi  \in H^1(0,T;H) \cap L^\infty (0,T;V) \cap  L^2(0,T;W), 
\end{equation}
\begin{equation} \label{regolaritaesistenza2reg}
\varphi    \in W^{1, \infty}(0,T;V^*) \cap H^1(0,T;V) \cap L^\infty (0,T;W), 
\end{equation}
\begin{equation} \label{regolaritaesistenza3reg}
\mu        \in L^{\infty}(0,T;V) \cap L^2(0,T;W).
\end{equation}
\end{Teo11-Regolare}

\paragraph{Remark.}
We fix $t\in (0,T)$ and integrate \eqref{iniziale2} over $\Omega$. We infer that
\begin{equation} \label{integraleconservato}
\int_{\Omega} \partial_t \varphi(t) - \int_{\Omega} \Delta \mu (t) = 0.
\end{equation}
Integrating by parts the second term of the left-hand side of \eqref{integraleconservato}, we obtain that
\begin{equation} \label{utileinregolaritaitti}
\frac{d}{dt} \int_{\Omega} \varphi(t) = 0.
\end{equation}
Consequently we conclude that
\begin{equation} \label{integraleconservato11}
m(\varphi(t)) = \frac{1}{|\Omega|} \int_{\Omega} \varphi(t) = \frac{1}{|\Omega|} \int_{\Omega} \varphi_0 = m (\varphi_0) =: m_0
\quad \textrm{for every $t \in (0,T)$.}
\end{equation}

\paragraph{Change of variables.}
In the following it we will be useful to consider the equivalent modified form of the initial Problem ($P$) (see \eqref{iniziale1}--\eqref{iniziale7}). We make a change of variables and set 
\begin{equation} \label{vivcinissima}
\kappa = a\chi + b \varphi - \kappa^* , \quad \quad \kappa_0 = a\chi_0 + b \varphi_0 - \kappa^*. 
\end{equation} 
Due to \eqref{vivcinissima}, from \eqref{iniziale1}--\eqref{iniziale7} we obtain the modified problem ($\widetilde{P}$):
\begin{equation} \label{iniziale1-uni1}
\partial_t(\kappa + (a\ell -b)\varphi) - \Delta \kappa + b \Delta \varphi - \Delta \kappa^* + a \zeta = a f \ \ \textrm{a.e. in $Q$,} 
\end{equation}
\begin{equation} \label{iniziale2-uni1}
\partial_t \varphi - \Delta \mu= 0  \ \ \textrm{a.e. in $Q$,}
\end{equation}
\begin{equation} \label{iniziale3-uni1}
\mu= - \nu \Delta \varphi + \xi + \pi(\varphi) - \frac{\gamma}{a} (\kappa - b \varphi + \kappa^*) \ \ \textrm{a.e. in $Q$,} 
\end{equation}
\begin{equation} \label{iniziale4-uni1}
\zeta(t) \in A(\kappa(t)) \ \textrm{for a.e. $t \in (0,T)$,}  
\end{equation}
\begin{equation} \label{iniziale5-uni1}
\xi \in \beta(\varphi) \ \textrm{a.e. in $Q$,}
\end{equation}
\begin{equation} \label{iniziale6-uni1}
\partial_{\bf{n}} \kappa = \partial_{\bf{n}} \varphi = \partial_{\bf{n}} \mu = 0 \ \ \textrm{on $\Sigma$,}
\end{equation}
\begin{equation} \label{iniziale7-uni1}
\kappa(0) = \kappa_0, \quad \quad \varphi(0) = \varphi_0 \quad \textrm{in $\Omega$}.
\end{equation}

\begin{Teo2-Unicita-Dipendenzacontinua}[Uniqueness and continuous dependence] \label{Teorema-unicita} 
Assume \eqref{beta}--\eqref{pi}, \\ \eqref{A1}--\eqref{stimaA} and \eqref{parametri}--\eqref{0}. 
If $a, \ b > 0$ and $a \ell = b$, 
then the solution $(\kappa,\varphi,\mu)$ of problem $(\widetilde{P})$ (see \eqref{iniziale1-uni1}--\eqref{iniziale7-uni1})
is unique. Moreover, we assume that $f_i$, $\kappa^{*}_i$, $\kappa_{0_i}$, $\varphi_{0_i}$, $i=1,2$, 
are given as in \eqref{f}--\eqref{0} and $(\kappa_i, \varphi_i,\mu_i)$, $i=1,2$,  are the corresponding solutions.
If 
\begin{equation} \label{masseuguali}
m(\varphi_{0_1}) = m(\varphi_{0_2}),
\end{equation} 
then the estimate
\begin{equation} \nonumber
\|\kappa_1 - \kappa_2\|_{L^{\infty}(0,T;H) \cap L^2(0,T;V)   } + \| \varphi_1 - \varphi_2\|_{L^{\infty}(0,T;V^*) \cap L^2(0,T;V)  }
\end{equation} 
\begin{equation} \label{riferimentodipendenzacontinua}
\leq  c \big( \| \varphi_{0_1} - \varphi_{0_2} \|_{V^*} + \|\kappa_{0_1} - \kappa_{0_2}\|_H + \|f_1 - f_2\|_{L^2(0,T; H)} +  \| \kappa^*_1- \kappa^*_2\|_W \big)   
\end{equation} 
holds true for some constant $c$ that depends only on $\Omega$, $T$
and the structure \eqref{beta}--\eqref{pi}, \eqref{A1}--\eqref{stimaA} and \eqref{parametri}--\eqref{0} of the system.
\end{Teo2-Unicita-Dipendenzacontinua}

\begin{Teo3-Sliding}[Sliding mode control] \label{Teo3-Sliding-mode}
Assume \eqref{beta}--\eqref{pi}, \eqref{A1}--\eqref{stimaA}, \eqref{parametri},
$a=1$, $ b = \ell$ and
\begin{equation} 
f \in L^{\infty}(0,T, H),
\end{equation}
\begin{equation}
\kappa^{*} \in W, \quad  \chi_0 \in V, \quad  \varphi_0 \in W,
\quad \beta^0(\varphi_0) \in H, \quad m_0 \in \interna(D(\beta)).
\end{equation}
We consider $A = \rho\Sign$, where $\rho$ is a positive coefficient, $\Sign$ is defined as in \eqref{Sign}
and $\sigma$ is an element of the range of $\Sign$, i.e.,
\begin{equation} \label{riferimentodefrho}
\sigma(t) \in \Sign(\chi(t)+\ell\varphi(t) - \kappa^*) \ \textrm{for a.e. $t \in (0,T)$,}
\end{equation}
Then, for some $\rho^* > 0$ and for every $\rho >\rho^*$, there exists a solution $(\chi,\varphi,\mu)$
to Problem $(P)$ (see \eqref{iniziale1}--\eqref{iniziale7}) and a time $T^*$ such that, for every $t \in [T^*,T]$
\begin{equation} \label{finedellarticolo3}
\chi(t) + \ell \varphi(t) = \kappa^* \quad \quad \textrm{a.e. in $\Omega$.}
\end{equation}
\end{Teo3-Sliding}

\setcounter{equation}{0}
\section{Existence - The approximating problem ($P_{\varepsilon}$)}
The following three sections are devoted to the proof of the existence Theorem \ref{Teorema-esistenza}.

Let us stress that, from now on, the symbol $c$ stands for different positive constants which depend only on $|\Omega|$, on the final time $T$,
the shape of the nonlinearities and on the constants and the norms of the functions involved in the assumptions of our statements.

\paragraph{Yosida regularization of $A$.} 
We consider the Yosida regularization of $A$. For $\varepsilon>0$ we define
\begin{equation} \label{richiamoA1}
A_{\varepsilon}: H \longrightarrow  H, \quad \quad
A_{\varepsilon} = \frac{I - (I + \varepsilon A)^{-1}}{\varepsilon},
\end{equation}
where $I$ denotes the identity operator.
Note that $A_{\varepsilon}$ is Lipschitz continuous and maximal monotone,
with Lipschitz constant $1/ \varepsilon$,
and satisfies the following properties. Denoting by $J_{\varepsilon}= (I + \varepsilon A)^{-1}$ the resolvent operator, for all $\delta > 0$ we have that
	\begin{equation} \label{AJe}
	A_{\varepsilon}\kappa \in A(J_{\varepsilon}\kappa),
	\end{equation}
	\begin{equation}
	(A_{\varepsilon})_\delta = A_{\varepsilon + \delta},
	\end{equation}
	\begin{equation}
	\| A_{\varepsilon}\kappa \|_H \leq \| A^0 \kappa \|_H,
	\end{equation}
	\begin{equation}
	\lim_{\varepsilon \rightarrow 0} \| A_{\varepsilon}\kappa \|_H = \| A^0 \kappa \|_H,
	\end{equation}
where $A^0 \kappa$ is the element	of the range of $A \kappa$ having minimum norm. 

\paragraph{Remark.}
We point out a key property of $A_{\varepsilon}$, 
which is a consequence of \eqref{stimaA}:
\begin{equation} \label{Ae}
\| v \|_H \leq C_A (1 + \| \kappa \|_H ) \quad \quad \textrm{for all} \ \kappa \in H, \ v \in A_{\varepsilon} \kappa.
\end{equation}
Indeed notice that $0 \in A(0)$ and $0 \in I(0)$: consequently, for every $\varepsilon >0$, $0 \in (I + \varepsilon A)(0)$. This fact implies that
$ J_\varepsilon (0)=0$.
We also recall that $A$ is a maximal monotone operator, hence $J_\varepsilon$ is a contraction. Then, from \eqref{stimaA} and \eqref{AJe}, it follows that
\begin{equation} \nonumber
\begin{array}{ll}
\| A_{\varepsilon}\kappa \|_H & \leq C_A (\| J_{\varepsilon}\kappa\|_H  + 1)\\
& \leq C_A (\| J_{\varepsilon}\kappa - J_{\varepsilon}0 \|_H + \| J_{\varepsilon}0 \|_H + 1)\\
& \leq C_A (\| \kappa \|_H + 1).
\end{array}													
\end{equation}

\paragraph{Yosida regularization of $\Sign$.} 
Let us introduce the operator $\Signe: \ H \rightarrow H$ as the Yosida regularization at level $\varepsilon>0$
of the operator $\Sign$.
We observe that $\Signe(v)$ is the gradient at $v $ of the $C^1$ functional  $\| \cdot \|_{H, \varepsilon}$
defined as
\begin{equation} \label{segno1}
\| v \|_{H, \varepsilon} := \min_{w \in H}{\{ \frac{1}{2\varepsilon} \| w-v \|^2_H + \| w \|_H     \}} 
= \int_0^{\| v \|_H} \min{\{ s/ \varepsilon, 1 \}} \ ds
\quad \textrm{for every $v \in H$}.
\end{equation} 
We also recall that
\begin{equation}  \label{segno2}
\Signe (v) = \frac{v}{\max{\{ \varepsilon, \|v \|_H \}}} \ \ \textrm{for every $v \in H$}.
\end{equation} 

\paragraph{Moreau-Yosida regularization of $\beta$ and $\widehat{\beta}$.} 
We introduce the Yosida regularization of $\beta$. For every $\varepsilon>0$ we define
\begin{equation} \label{richiamobeta1}
\beta_{\varepsilon}: \mathbb{R} \longrightarrow \mathbb{R},
\quad \quad
\beta_{\varepsilon} = \frac{I - (I + \varepsilon \beta)^{-1}}{\varepsilon}. 
\end{equation}
We remark that $\beta_{\varepsilon}$ is Lipschitz continuous with Lipschitz constant $1/ \varepsilon$ and satisfies the following properties. Denoting by $R_{\varepsilon}= (I + \varepsilon \beta)^{-1}$ the resolvent operator, for all $\delta > 0$ and for every $\varphi \in D(\beta)$ we have that
	\begin{equation}
	\beta_{\varepsilon}(\varphi) \in \beta(R_{\varepsilon}\varphi),
  \end{equation}
		\begin{equation}
	(\beta_{\varepsilon})_\delta = \beta_{\varepsilon + \delta}, 
	\end{equation}
	\begin{equation}
	|\beta_{\varepsilon}(\varphi)| \leq |\beta^0(\varphi)|,
	\end{equation}
	\begin{equation}
 	\lim_{\varepsilon \rightarrow 0} \beta_{\varepsilon}(\varphi) =  \beta^0(\varphi), 
	\end{equation}
where $\beta^0(\varphi)$ is the element	of the range of $\beta$ having minimum modulus.
For $\varepsilon>0$, we also introduce $\widehat{\beta}_{\varepsilon}: \mathbb{R} \rightarrow [0, +\infty]$ 
as the standard Moreau-Yosida regularization of $\widehat{\beta}$
\begin{equation}
\widehat{\beta}_{\varepsilon} := \min_{y \in \mathbb{R}} \bigg\{  \widehat{\beta}(x) + \frac{1}{2 \varepsilon} |x - y| \bigg\}
\end{equation}
and we recall that, for every $\varphi  \in D(\widehat{\beta})$,
	\begin{equation} \label{questaservepoi}
	\widehat{\beta}_{\varepsilon}(\varphi) \leq \widehat{\beta}(\varphi).
	\end{equation}
Moreover, $\beta_{\varepsilon}$ is the Frechet derivative of $\widehat{\beta}_{\varepsilon}$. Then, for every $\varphi_1, \varphi_2  \in D(\widehat{\beta})$, we have that
\begin{equation} \label{questaservepoimolto}
\widehat{\beta}_{\varepsilon}(\varphi_2) = \widehat{\beta}_{\varepsilon}(\varphi_1) + \int_{\varphi_1}^{\varphi_2} \beta_{\varepsilon}(s) \ ds.
\end{equation}

\paragraph{Regularization of the initial data.}
We denote  by $\chi_{0\varepsilon}$ and 
$\varphi_{0\varepsilon}$  the regularization of the initial data $\chi_0$ and $\varphi_0$, respectively,
obtained solving the following elliptic problems:
\begin{equation} \label{sitemadatiinizialepsilon1}
\left\{ \begin{array}{ll}
\chi_{0\varepsilon} - \varepsilon \Delta \chi_{0\varepsilon} = \chi_0							& \textrm{in $\Omega$}, \\
\partial_{\bf{n}}\chi_{0\varepsilon} = 0																												& \textrm{on $ \Gamma  $.}
\end{array}
\right.
\end{equation}
\begin{equation} \label{sitemadatiinizialepsilon2}
\left\{ \begin{array}{ll}
\varphi_{0\varepsilon} - \varepsilon \Delta \varphi_{0\varepsilon} = \varphi_0							& \textrm{in $\Omega$}, \\
\partial_{\bf{n}}\varphi_{0\varepsilon} = 0																												& \textrm{on $ \Gamma  $.}
\end{array}
\right.
\end{equation}
Since $\chi_0 \in H$ and $\varphi_0 \in V$, by elliptic regularity we infer that
$\chi_{0\varepsilon} \in W$ and $\varphi_{0\varepsilon} \in W \cap H^3(\Omega)$.
Moreover, integrating over $\Omega$ the first equation of \eqref{sitemadatiinizialepsilon2},  we obtain that 
\begin{equation} \label{eluiilveroriferimentoresatto}
m_0= \frac{1}{| \Omega |} \int_{\Omega } \varphi_0 = \frac{1}{| \Omega |} \int_{\Omega } \varphi_{0  \varepsilon} =: m_{0\varepsilon}.
\end{equation} 
From \eqref{0} and \eqref{integraleconservato11} it immediately follows that 
$m_{0\varepsilon} \in \interna (D(\beta))$. 
Since $\beta$ is maximal monotone, testing the first equation of \eqref{sitemadatiinizialepsilon2} 
by $\beta_{\varepsilon}(\varphi_{0\varepsilon})$ and
integrating over $\Omega$, we have that
\begin{equation} \label{b1111}
\int_{\Omega} (\varphi_{0\varepsilon} - \varphi_0  ) \beta_{\varepsilon}(\varphi_{0\varepsilon})
= - \varepsilon \int_{\Omega} | \nabla \varphi_{0\varepsilon} |^2  \beta'_{\varepsilon}(\varphi_{0\varepsilon}) \leq 0.
\end{equation}
Recalling that $\beta_{\varepsilon}$ is the subdifferential of $\widehat{\beta}_{\varepsilon}$, from \eqref{b1111}
we infer that
\begin{equation} \label{stimanecessary}
\int_{\Omega} \widehat{\beta}_{\varepsilon} (\varphi_{0\varepsilon}) - \int_{\Omega} \widehat{\beta}_{\varepsilon} (\varphi_0) 
\leq \int_{\Omega} (\varphi_{0\varepsilon} - \varphi_0) \beta_{\varepsilon}(\varphi_{0\varepsilon}) \leq 0.
\end{equation}
Consequently,  due to \eqref{0},  \eqref{questaservepoi}, \eqref{stimanecessary} and the definition of $\widehat{\beta}_{\varepsilon}$,
we conclude that
\begin{equation}
0 \leq \int_{\Omega} \widehat{\beta}_{\varepsilon} (\varphi_{0\varepsilon}) 
\leq \int_{\Omega} \widehat{\beta}_{\varepsilon} (\varphi_0) 
\leq \int_{\Omega} \widehat{\beta} (\varphi_0) < + \infty,
\end{equation}
whence there exists a positive constant $c$, independent of $\varepsilon$,
such that $\| \widehat{\beta}(\varphi_{0\varepsilon}) \|_{ L^1(\Omega)} \leq c$.
Now, we test \eqref{sitemadatiinizialepsilon1}
by $\chi_{0\varepsilon}$ and integrate over $\Omega$. We obtain that
\begin{equation} \label{4444-123}
\int_{\Omega} |\chi_{0\varepsilon}|^2 + \varepsilon \int_{\Omega} |\nabla \chi_{0\varepsilon}|^2 
= \int_{\Omega} \chi_0 \chi_{0\varepsilon} \leq \frac{1}{2} \int_{\Omega} |\chi_0|^2 +
 \frac{1}{2}  \int_{\Omega} | \chi_{0\varepsilon}|^2.
\end{equation}
Since $\chi_0 \in H$, from \eqref{4444-123} it immediately follows that
$ \varepsilon \chi_{0\varepsilon} \longrightarrow 0$ in $V$ as $\varepsilon \searrow 0$.
Besides, there exists a positive constant $c$,
independent of $\varepsilon$, such that $\| \chi_{0\varepsilon} \|_H \leq c$. 
Then, testing the first equation of the system \eqref{sitemadatiinizialepsilon1} by an arbitrary function $v \in V$
and passing to the limit as $\varepsilon \searrow 0$, we obtain that
\begin{equation} \label{ccc1}
\lim_{\varepsilon \searrow 0} \Bigg( \int_{\Omega} \chi_{0\varepsilon} v + \varepsilon \int_{\Omega} \nabla \chi_{0\varepsilon} \cdot \nabla v - \int_{\Omega} \chi_0		v \Bigg)	= 0 	\quad \quad			 \textrm{for all $v \in V$}, 
\end{equation}
whence $	\chi_{0\varepsilon} \rightharpoonup \chi_0 $ in $H$.
Moreover, from \eqref{4444-123} and \eqref{ccc1} we infer that
\begin{equation} \label{convnormteta}
\int_{\Omega} |\chi_{0}|^2 
\leq \liminf_{\varepsilon \searrow 0} \int_{\Omega} |\chi_{0\varepsilon}|^2 
\leq \limsup_{\varepsilon \searrow 0} \int_{\Omega} |\chi_{0\varepsilon}|^2 
\leq \int_{\Omega} |\chi_0|^2.
\end{equation}
Thanks to \eqref{convnormteta}, $ \|\chi_{0\varepsilon} \|_H \longrightarrow \|\chi_0 \|_H$
and this ensures, due to the weak convergence already proved, that $ \chi_{0\varepsilon} \longrightarrow \chi_0 $ in $H$.

With a similar technique, testing \eqref{sitemadatiinizialepsilon2}
by $\varphi_{0\varepsilon}$ and integrating over $\Omega$, we obtain that $ \varphi_{0\varepsilon} \longrightarrow \varphi_0 $ in $H$.
Now, we test \eqref{sitemadatiinizialepsilon2}
by $- \Delta \varphi_{0\varepsilon}$ and integrate over $\Omega$. We obtain that
\begin{equation} \label{4444-123fd}
\int_{\Omega} |\nabla \varphi_{0\varepsilon}|^2 + \varepsilon \int_{\Omega} |\Delta \varphi_{0\varepsilon}|^2 
= \int_{\Omega} \nabla \varphi_0 \cdot \nabla \varphi_{0\varepsilon} \leq \frac{1}{2} \int_{\Omega} | \nabla \varphi_0|^2 +
 \frac{1}{2}  \int_{\Omega} | \nabla \varphi_{0\varepsilon}|^2.
\end{equation}
Since $\varphi_0 \in V$, from \eqref{4444-123fd} it immediately follows that
$ \varepsilon \varphi_{0\varepsilon} \longrightarrow 0$ in $W$ as $\varepsilon \searrow 0$.
Furthermore, there exists a positive constant $c$,
independent of $\varepsilon$, such that $\| \nabla \varphi_{0\varepsilon} \|_H \leq c$.
Recalling that  $\| \varphi_{0\varepsilon} \|_H \leq c$, we conclude that $\| \varphi_{0\varepsilon} \|_V \leq c$.
Then, testing the the first equation of the system \eqref{sitemadatiinizialepsilon2} by $- \Delta w$, where $w$
is an arbitrary function in $W$, and passing to the limit as $\varepsilon \searrow 0$, we obtain
\begin{equation} \label{cstr2}
\lim_{\varepsilon \searrow 0} \Bigg( \int_{\Omega} \nabla \varphi_{0\varepsilon} \cdot \nabla w + \varepsilon \int_{\Omega} \Delta \varphi_{0\varepsilon} \cdot \Delta w - \int_{\Omega}  \nabla \varphi_0	\cdot \nabla	w \Bigg)	= 0 	\quad \quad			 \textrm{for all $w \in W$}, 
\end{equation}
whence $	\varphi_{0\varepsilon} \rightharpoonup \varphi_0 $ in $V$.
Moreover, from \eqref{4444-123fd}--\eqref{cstr2} we infer that
\begin{equation} \label{convnormtetafi}
\int_{\Omega} |\nabla \varphi_0|^2
\leq \liminf_{\varepsilon \searrow 0} \int_{\Omega} | \nabla \varphi_{0\varepsilon}|^2 
\leq \limsup_{\varepsilon \searrow 0} \int_{\Omega} | \nabla \varphi_{0\varepsilon}|^2 
\leq \int_{\Omega} |\nabla \varphi_0|^2.
\end{equation}
Thanks to \eqref{convnormtetafi}, $ \| \nabla \varphi_{0\varepsilon} \|_H \longrightarrow \|\nabla \varphi_0 \|_H$
and this ensures, due to the weak convergence already proved, that $ \varphi_{0\varepsilon} \longrightarrow \varphi_0 $ in $V$.
Now, let us summarize the main properties of $\chi_{0\varepsilon}$ and $\varphi_{0\varepsilon}$.
For every $\varepsilon \in (0,1)$ we have that
\begin{equation} \label{0epsilon}
\chi_{0\varepsilon} \in W,  \quad 
\varphi_{0\varepsilon} \in W \cap H^3(\Omega),  \quad 
m_{0\varepsilon} \in \interna (D(\beta)),  \quad
\| \widehat{\beta}(\varphi_{0\varepsilon}) \|_{ L^1(\Omega)} \leq c,
\end{equation}
\begin{equation} \label{convergenzadatiinizialiepsilon1} 
\lim_{\varepsilon \searrow 0} \| \chi_0  - \chi_{0\varepsilon} \|_H =0,
\quad \quad
\lim_{\varepsilon \searrow 0} \| \varphi_0  - \varphi_{0\varepsilon} \|_V =0,
\end{equation}
\begin{equation} \label{tecnicadistima2}
-\nu \Delta \varphi_{0\varepsilon} + \beta_{\varepsilon}(\varphi_{0\varepsilon}) 
+ \pi(\varphi_{0\varepsilon}) - \gamma \chi_{0\varepsilon}  \in V .
\end{equation}

\paragraph{Regularization of $f$.}
We denote by $f_{\varepsilon}$ the regularization of $f$, constructed in such a way that
\begin{equation} \label{convergenzadif} 
f_{\varepsilon} \in C^1([0,T]; H) \textrm{ for all $\varepsilon >0$}, \quad \quad
\lim_{\varepsilon \searrow 0} \| f_{\varepsilon} - f \|_{L^2(0,T; H)} = 0.
\end{equation}
For example, we can consider $f_{\varepsilon}$ as the solution the following system:
\begin{equation}
\left\{ \begin{array}{ll}
-\varepsilon f''_{\varepsilon}(t) + f_{\varepsilon}(t)  = f(t), \quad t \in (0,T),  \\
f_{\varepsilon}(0) = f_{\varepsilon}(T) = 0. \\
\end{array}
\right.
\end{equation}

\paragraph{Approximating problem ($P_{\varepsilon}$).}
We look for a triplet $(\chi_{ \varepsilon},\varphi_{ \varepsilon},\mu_{ \varepsilon})$ satisfying at least the regularity requirements 
\begin{equation} \label{regolaritaesistenza-epsilon1}
\chi_{\varepsilon}  \in H^1(0,T;H) \cap L^\infty (0,T;V) \cap  L^2(0,T;W), 
\end{equation}
\begin{equation} \label{regolaritaesistenza-epsilon2}
\varphi_{\varepsilon}    \in W^{1, \infty}(0,T;V^*) \cap H^1(0,T;V) \cap L^\infty (0,T;W), 
\end{equation}
\begin{equation} \label{regolaritaesistenza-epsilon3}
\mu_{\varepsilon}        \in L^{\infty}(0,T;V) \cap L^2(0,T;W),
\end{equation}
and solving the approximating problem ($P_{\varepsilon}$):
\begin{equation} \label{iniziale1-epsilon}
\partial_t(\chi_{\varepsilon} + \ell \varphi_{\varepsilon}) - \Delta \chi_{\varepsilon} + \zeta_{\varepsilon} = f_{\varepsilon} \ \ \textrm{a.e. in $Q$,} 
\end{equation}
\begin{equation} \label{iniziale2-epsilon}
\partial_t \varphi_{\varepsilon} - \Delta \mu_{\varepsilon}= 0  \ \ \textrm{a.e. in $Q$,}
\end{equation}
\begin{equation} \label{iniziale3-epsilon}
\mu_{\varepsilon}= - \nu \Delta \varphi_{\varepsilon} + \xi_{\varepsilon} + \pi(\varphi_{\varepsilon}) - \gamma \chi_{\varepsilon} \ \ \textrm{a.e. in $Q$,} 
\end{equation}
\begin{equation} \label{iniziale4-epsilon}
\zeta_{\varepsilon}(t) \in A_{\varepsilon}(a \chi_{\varepsilon}(t) + b \varphi_{\varepsilon}(t) - \kappa^*) \ \textrm{for a.e. $t \in (0,T)$,}  
\end{equation}
\begin{equation} \label{iniziale5-epsilon}
\xi_{\varepsilon} \in \beta_{\varepsilon}(\varphi_{\varepsilon}) \ \textrm{a.e. in $Q$,}
\end{equation}
\begin{equation} \label{iniziale6-epsilon}
\partial_{\bf{n}} \chi_{\varepsilon} = \partial_{\bf{n}} \varphi_{\varepsilon} = \partial_{\bf{n}} \mu_{\varepsilon} = 0 \ \ \textrm{on $\Sigma$,}
\end{equation}
\begin{equation} \label{iniziale7-epsilon}
\chi_{\varepsilon}(0) = \chi_{0\varepsilon}, \quad \quad \varphi_{\varepsilon}(0) = \varphi_{0\varepsilon} \ \ \textrm{in $\Omega$,}
\end{equation}
where $\beta_{\varepsilon}$ and $A_{\varepsilon}$ are the Yosida regularizations of $\beta$ and $A$
defined in \eqref{richiamoA1} and \eqref{richiamobeta1}.
We notice that the homogeneous Neumann boundary conditions
are already contained in the conditions \eqref{regolaritaesistenza-epsilon1}--\eqref{regolaritaesistenza-epsilon3}
due to the definition of $W$ (see \eqref{W}).

We observe that, for almost every $t \in (0, T)$, we can re-write the approximating problem ($P_{\varepsilon}$) 
in the following way:
\begin{equation}\label{riferimentostima1}
\langle \partial_t(\chi_{\varepsilon} + \ell \varphi_{\varepsilon}) (t), z \rangle_{V^*,V}
+ \int_{\Omega} \nabla \chi_{\varepsilon}(t) \cdot \nabla z
+ \langle \zeta_{\varepsilon}(t), z \rangle_{V^*,V} = \langle f_{\varepsilon}(t), z \rangle_{V^*,V} \ \ \textrm{for all $z \in V$,} 
\end{equation}
\begin{equation} \label{riferimentostima2}
\langle  \partial_t \varphi_{\varepsilon}(t), z \rangle_{V^*,V} + \int_{\Omega} \nabla \mu_{\varepsilon}(t) \cdot \nabla z = 0  \ \ \textrm{for all $z \in V$,} 
\end{equation}
\begin{equation} \label{riferimentostima3}
\mu_{\varepsilon}(t)= - \nu \Delta \varphi_{\varepsilon}(t) + \xi_{\varepsilon}(t) + \pi(\varphi_{\varepsilon}(t)) - \gamma \chi_{\varepsilon}(t) \ \ \textrm{in $H$,} 
\end{equation}
\begin{equation} \label{iniziale422}
\zeta_{\varepsilon}(t) \in A_{\varepsilon}(a \chi_{\varepsilon}(t) + b \varphi_{\varepsilon}(t) - \kappa^*),  
\end{equation}
\begin{equation} \label{iniziale522}
\xi \in \beta_{\varepsilon}(\varphi_{\varepsilon}) \ \textrm{a.e. in $Q$,}
\end{equation}
\begin{equation} \label{iniziale622}
\partial_{\bf{n}} \varphi_{\varepsilon} = 0 \quad \textrm{a.e. on $\Sigma$,}
\end{equation}
\begin{equation} \label{iniziale722}
\chi_{\varepsilon}(0) = \chi_{0\varepsilon}, \quad \quad \varphi_{\varepsilon}(0) = \varphi_{0\varepsilon} 
\quad \textrm{in $\Omega$.}
\end{equation}
Since $m_{0\varepsilon}=m_0$, 
recalling the definition of $\mathcal{N}$ (see \eqref{operatoreN1}--\eqref{operatoreN3}), we have that $\partial_t \varphi_{\varepsilon} (t) \in D(\mathcal{N})$. Hence, \eqref{riferimentostima2} can be written as
\begin{equation} \label{riferimentostima4}
\mathcal{N}\partial_t \varphi_{\varepsilon} (t) = m(\mu_{\varepsilon}(t)) - \mu_{\varepsilon}(t) \ \ \textrm{in $V$,}
\end{equation}
and this and \eqref{riferimentostima2} entail
\begin{equation} \label{riferimentostima5}
m(\mu_{\varepsilon}(t)) - \mathcal{N}\partial_t \varphi_{\varepsilon} (t)
=  - \nu \Delta \varphi_{\varepsilon}(t) + \xi_{\varepsilon}(t) + \pi(\varphi_{\varepsilon}(t)) - \gamma \chi_{\varepsilon}(t) \ \ \textrm{in $H$.} 
\end{equation}

\setcounter{equation}{0}
\section{Existence - Global a priori estimates}
In this se ction, we  will de duce  some  a priori e stimate s 
infe rre d from \eqref{riferimentostima1}--\eqref{riferimentostima5}.

In the  re mainde r of the  pape r we  ofte n owe  to the  H\"olde r ine quality and to the  e le me ntary Young ine qualitie s in pe rforming our a priori e stimate s.
For e very $x,y >0$, $\alpha \in (0,1)$ and $\delta > 0$ there hold
\begin{equation} \label{dis1}
xy \leq \alpha x^{\frac{1}{\alpha}} + (1 - \alpha)y^{\frac{1}{1- \alpha}},
\end{equation}
\begin{equation} \label{dis2}
xy \leq \delta x^2 + \frac{1}{4\delta}y^2.
\end{equation}
More ove r, we  also use  the  ine quality de duce d from the  compactne ss 
of the  e mbedding $V \subset H \subset V^*$ (see \cite[Lemma 8, p. 84]{Sim87}): for all $\delta > 0$ 
the re  e xists a constant $K > 0$ such that 
\begin{equation} \label{stimasimon1}
\| z \|_H 
\leq \delta \| z \|_V + K \| z \|_{V^*} 
\quad \quad \textrm{for all $z \in H$.}
\end{equation}
In the  following, the  symbol $c$ stands for diffe re nt positive  constants which de pe nd only on $|\Omega|$, on the final time $T$,
on the  shape  of the  nonline aritie s and on the  constants and the  norms of the  functions involve d in the  assumptions of our state me nts.

\paragraph{First a priori estimate.}
According to \eqref{eluiilveroriferimentoresatto}, $m(\partial_t \varphi_{\varepsilon})=0$. Consequently,
$\partial_t \varphi_{\varepsilon} \in D(\mathcal{N})$ and we can test \eqref{riferimentostima2} 
by $\mathcal{N}\partial_t\varphi_{\varepsilon}$. Integrating over $(0,t)$, $t \in (0,T]$, we obtain that
\begin{equation} \label{preambolo}
\int_0^t \| \partial_t\varphi_{\varepsilon}(s) \|^2_{V^*} \ ds + \int_{Q_t} \nabla \mu_{\varepsilon} \cdot \nabla \mathcal{N} \partial_t\varphi_{\varepsilon} 
= \int_0^t \| \partial_t\varphi_{\varepsilon} \|^2_{V^*} + \int_{Q_t} \mu_{\varepsilon} \partial_t\varphi_{\varepsilon} = 0.
\end{equation} 
Recalling that
\begin{equation} \label{aggiunta}
\nu \int_{Q_t} \varphi_{\varepsilon} \partial_t \varphi_{\varepsilon} = \frac{\nu}{2} \int_{\Omega} |\varphi_{\varepsilon} (t)|^2 - \frac{\nu}{2} \int_{\Omega} |\varphi_{0\varepsilon}|^2,
\end{equation} 
we combine \eqref{riferimentostima1} tested by $\frac{\gamma}{\ell}\chi_{\varepsilon}$, \eqref{preambolo} and \eqref{aggiunta}.
Then we subtract \eqref{riferimentostima3} tested by $\partial_t\varphi_{\varepsilon}$ and integrate over $(0,t)$. 
We have that
\begin{equation} \nonumber
\frac{\gamma}{2\ell} \int_{\Omega} |\chi_{\varepsilon}(t)|^2 
+ \frac{\gamma}{\ell} \int_{Q_t} |\nabla \chi_{\varepsilon}|^2 
+ \int_0^t \| \partial_t\varphi_{\varepsilon} (s)\|^2_{V^*} \ ds
+ \frac{\nu}{2} \| \varphi_{\varepsilon}(t) \|^2_V 
+ \int_{\Omega} \widehat{\beta_{\varepsilon}}(\varphi_{\varepsilon}(t))
\end{equation} 
\begin{equation}\label{prima1}
= \frac{\gamma}{2 \ell} \|\chi_{0\varepsilon}\|_H^2 + \frac{\nu}{2} \| \varphi_{0\varepsilon} \|^2_V + \int_{\Omega} \widehat{\beta_{\varepsilon}}(\varphi_{0\varepsilon})
+ \frac{\gamma}{\ell} \int_{Q_t} f_{\varepsilon} \chi_{\varepsilon} 
- \frac{\gamma}{\ell}\int_{Q_t} \zeta_{\varepsilon} \chi_{\varepsilon}
+ \int_{Q_t} (\nu \varphi_{\varepsilon} - \pi(\varphi_{\varepsilon}))\partial_t\varphi_{\varepsilon}.
\end{equation} 
As $\pi$ is a Lipschitz continuous function with Lipschitz constant $C_{\pi} = \| \pi' \|_{\infty}$, we obtain that 
\begin{equation} \label{stimainizailesupigreco}
|\pi(\varphi_{\varepsilon}(s))| 
\leq |\pi(\varphi_{\varepsilon}(s)) - \pi(0)| + |\pi(0)| 
\leq C_{\pi} |\varphi_{\varepsilon}(s)| + |\pi(0)|. 
\end{equation} 
Consequently, thanks to \eqref{stimainizailesupigreco}, we infer that
\begin{equation} \nonumber
\| \nu \varphi_{\varepsilon} (s) - \pi(\varphi_{\varepsilon}(s)) \|^2_V 
= \int_{\Omega} | \nu \varphi_{\varepsilon} (s) - \pi(\varphi_{\varepsilon}(s)) |^2 + \int_{\Omega} | \nu \nabla \varphi_{\varepsilon} (s) - \pi'(\varphi_{\varepsilon}(s)) \nabla \varphi_{\varepsilon}(s) |^2
\end{equation}
\begin{equation} \nonumber
\leq 2 \int_{\Omega} \bigg( \nu^2 |\varphi_{\varepsilon} (s)|^2 + |\pi(\varphi_{\varepsilon}(s)) |^2 \bigg) 
   + 2 \int_{\Omega} \bigg( \nu^2 |\nabla \varphi_{\varepsilon} (s)|^2 +  \| \pi' \|^2_{\infty} |\nabla \varphi_{\varepsilon}(s) |^2 \bigg)
\end{equation}
\begin{equation} \nonumber
\leq 2\nu^2 \int_{\Omega} | \varphi_{\varepsilon} (s)|^2 + 4 C^2_{\pi} \int_{\Omega}  |\varphi_{\varepsilon}(s)|^2 
+ 4 |\Omega| |\pi(0)|^2 + 2 \nu^2 \int_{\Omega} |\nabla \varphi_{\varepsilon} (s)|^2 + 2 C_{\pi}^2 \int_{\Omega}  |\nabla \varphi_{\varepsilon}(s) |^2
\end{equation} 
\begin{equation} \label{stimautilesupi}
= ( 2\nu^2 + 4 C^2_{\pi} ) \int_{\Omega} | \varphi_{\varepsilon} (s)|^2  + (2 \nu^2 + 2 C_{\pi}^2) \int_{\Omega}  |\nabla \varphi_{\varepsilon}(s) |^2
+ 4 |\pi(0)|^2 |\Omega|
\leq c ( \| \varphi_{\varepsilon} (s) \|^2_V  + 1).
\end{equation} 
Due to \eqref{stimautilesupi}, we obtain that the last term on the right-hand side of \eqref{prima1}
is estimated as follows
\begin{eqnarray} \nonumber
\int_{Q_t} (\nu \varphi_{\varepsilon} - \pi(\varphi_{\varepsilon})) \partial_t\varphi_{\varepsilon}
& \leq&  \frac{1}{2} \int_0^t \| \partial_t\varphi_{\varepsilon} (s)\|^2_{V^*} \ ds + \frac{1}{2} \int_0^t \| \nu \varphi_{\varepsilon} (s) - \pi(\varphi_{\varepsilon}(s)) \|^2_V \ ds \\ \label{pigrecostima}
& \leq & \frac{1}{2} \int_0^t \| \partial_t\varphi_{\varepsilon} (s)\|^2_{V^*} \ ds + c \int_0^t ( \| \varphi_{\varepsilon} (s) \|^2_V  + 1) \ ds . 
\end{eqnarray} 
Due to the liear growth of $A_{\varepsilon}$ stated by \eqref{Ae}, we have that
\begin{equation} \nonumber
- \frac{\gamma}{\ell} \int_{Q_t} \zeta_{\varepsilon} \chi_{\varepsilon} 
\leq \frac{\gamma}{\ell}   \int_{Q_t} |\zeta_{\varepsilon}(s)| | \chi_{\varepsilon}(s)| \ ds 
\leq \frac{\gamma}{\ell}   \int_0^t \| \zeta_{\varepsilon}(s) \|^2_H \ ds + \frac{\gamma}{\ell}  \int_0^t \|\chi_{\varepsilon}(s)\|^2_H \ ds
\end{equation} 
\begin{equation} \nonumber
\leq \frac{\gamma}{\ell}   \int_0^t C_A^2 ( 1 + \| a\chi_{\varepsilon}(s) + b \varphi_{\varepsilon}(s) - \kappa^* \|_H)^2 \ ds + \frac{\gamma}{\ell}  \int_0^t \|\chi_{\varepsilon}(s)\|^2_H \ ds 
\end{equation} 
\begin{equation} \nonumber
\leq \frac{\gamma}{\ell}   \int_0^t 4C_A^2 ( 1 + |a|^2 \| \chi_{\varepsilon}(s) \|_H^2 +  |b|^2  \| \varphi_{\varepsilon}(s) \|_H^2 + \| \kappa^* \|_H^2) \ ds + \frac{\gamma}{\ell}  \int_0^t \|\chi_{\varepsilon}(s)\|^2_H \ ds 
\end{equation} 
\begin{equation} \nonumber
\leq \frac{\gamma}{\ell} 4C_A^2 T + \frac{\gamma}{\ell} 4C_A^2 |a|^2 \int_0^t \| \chi_{\varepsilon}(s) \|_H^2 \ ds + \frac{\gamma}{\ell} 
4C_A^2 |b|^2 \int_0^t \| \varphi_{\varepsilon}(s) \|_H^2 \ ds
\end{equation} 
\begin{equation} \nonumber
+ \frac{\gamma}{\ell} 4C_A^2 T \| \kappa^* \|_H^2 + \frac{\gamma}{\ell}  \int_0^t \|\chi_{\varepsilon}(s)\|^2_H \ ds 
\end{equation} 
\begin{equation} \label{Astima}
\leq c \Bigg( \int_0^t \| \chi_{\varepsilon}(s) \|_H^2 \ ds + \int_0^t \| \varphi_{\varepsilon}(s) \|_H^2 \ ds + 1 \Bigg).
\end{equation} 
Moreover, by applying \eqref{dis2} to the fourth term on the right-hand side of \eqref{prima1}, we have that
\begin{equation} \label{fstima}
\frac{\gamma}{\ell} \int_{Q_t} f_{\varepsilon} \chi_{\varepsilon} 
\leq \frac{\gamma}{\ell} \int_{Q_t} |f_{\varepsilon}|^2 + \frac{\gamma}{4\ell} \int_{Q_t} |\chi_{\varepsilon}|^2
= \frac{\gamma}{\ell} \int_{Q_t} |f_{\varepsilon}|^2 + \frac{\gamma}{4\ell} \int_0^t \|\chi_{\varepsilon}(s)\|_H^2 \ ds. 
\end{equation} 
We rearrange the right-hand side of \eqref{prima1} using \eqref{pigrecostima}--\eqref{fstima} and obtain that
\begin{equation} \nonumber
\frac{\gamma}{2\ell} \int_{\Omega} |\chi_{\varepsilon}(t)|^2 
+ \frac{\gamma}{\ell} \int_{Q_t} |\nabla \chi_{\varepsilon} |^2 
+ \frac{1}{2} \int_0^t \| \partial_t\varphi_{\varepsilon} (s)\|^2_{V^*} \ ds
+ \frac{\nu}{2} \| \varphi_{\varepsilon}(t) \|^2_V
+ \int_{\Omega} \widehat{\beta_{\varepsilon}}(\varphi_{\varepsilon}(t))
\end{equation} 
\begin{equation}\nonumber
\leq \frac{\gamma}{2 \ell} \|\chi_{0\varepsilon}\|_H^2 + \frac{\nu}{2} \| \varphi_{0\varepsilon} \|^2_V + \int_{\Omega} \widehat{\beta_{\varepsilon}}(\varphi_{0\varepsilon})
+ \frac{\gamma}{\ell} \int_0^t \|f_{\varepsilon} (s) \|_H^2 \ ds
\end{equation}
\begin{equation} \label{prima2}
+ c \Bigg( \int_0^t \| \varphi_{\varepsilon} (s) \|^2_V \ ds 
+ \int_0^t \| \chi_{\varepsilon}(s) \|_H^2 \ ds + 1 \Bigg) 
+ \frac{\gamma}{4\ell} \int_0^t \|\chi_{\varepsilon}(s)\|_H^2 \ ds.
\end{equation} 
Due to \eqref{0epsilon}--\eqref{convergenzadatiinizialiepsilon1},
the first three terms of the right-hand side of \eqref{prima2}
are bounded and similarly the fourth term, thanks to \eqref{convergenzadif}.
Then, applying the Gronwall lemma, we conclude that
there exists a positive constant $c$, independent of $\varepsilon$, such that
\begin{equation} \label{ultimissimalemma1}
\frac{\gamma}{2\ell} \int_{\Omega} |\chi_{\varepsilon}(t)|^2 + \int_{Q_t} |\nabla \chi_{\varepsilon}|^2 
+ \frac{1}{2} \int_0^t \| \partial_t\varphi_{\varepsilon} (s)\|^2_{V^*} \ ds
+ \frac{\nu}{2} \| \varphi_{\varepsilon}(t) \|^2_V
+ \int_{\Omega} \widehat{\beta_{\varepsilon}}(\varphi_{\varepsilon}(t)) \leq c,
\end{equation} 
whence  it immediately follows that
\begin{eqnarray}\label{primastimaaprioiri1}
\| \chi_{\varepsilon} \|_{L^{\infty}(0,T; H) \cap L^2(0,T; V)} & \leq  & c , \\ \label{primastimaaprioiri3}
\| \varphi_{\varepsilon} \|_{H^1(0,T; V^*) \cap L^{\infty}(0,T; V)}  & \leq  & c , \\ \label{primastimaaprioiri6}
\| \widehat{\beta_{\varepsilon}}(\varphi_{\varepsilon}) \|_{L^{\infty}(0,T; L^1(\Omega))}   & \leq  & c  . 
\end{eqnarray}
Due to \eqref{primastimaaprioiri1}--\eqref{primastimaaprioiri6}, 
by \eqref{Ae} we have that
\begin{equation} \label{primastimaaprioiri5}
\| \zeta_{\varepsilon} \|_{L^{\infty}(0,T; H)}  \leq c,
\end{equation}
and, consequently,
by comparison in \eqref{riferimentostima1} we infer that
\begin{equation} \label{primastimaaprioiri2}
\| \partial_t \chi_{\varepsilon} \|_{L^2(0,T; V^*)} \leq c  .
\end{equation}

\paragraph{Second a priori estimate.}
Recalling that $m_{0\varepsilon} = m_0$ due to \eqref{riferimentostima2},
we have that $\varphi_{\varepsilon}(s)-m_0 \in D(\mathcal{N})$
for every $s \in (0, T)$. We test \eqref{riferimentostima5} at time 
$s$ by $(\varphi_{\varepsilon}(s)-m_0) \in D(\mathcal{N})$ and we infer
that
\begin{equation} \nonumber
(\xi_{\varepsilon}(s), \varphi_{\varepsilon}(s)-m_0)_H=
- (\mathcal{N}\partial_t \varphi_{\varepsilon} (s), \varphi_{\varepsilon}(s)-m_0)_H
+ (m(\mu_{\varepsilon}(s)), \varphi_{\varepsilon}(s)-m_0)_H
\end{equation}
\begin{equation} \label{astrodiriferimento}
+ \nu (\Delta \varphi_{\varepsilon}(s), \varphi_{\varepsilon}(s)-m_0)_H
- (\pi(\varphi_{\varepsilon}(s)), \varphi_{\varepsilon}(s)-m_0)_H
+ \gamma (\chi_{\varepsilon}(s), \varphi_{\varepsilon}(s)-m_0)_H.
\end{equation}
We recall that there exists a positive constant $c$ such that 
$\|z\|_{V^*} \leq c \|z\|_H$ for all $z \in H$.
Consequently the first term of the right-hand side of \eqref{astrodiriferimento} is estimated as follows:
\begin{eqnarray} \nonumber
- (\mathcal{N}\partial_t \varphi_{\varepsilon} (s), \varphi_{\varepsilon}(s)-m_0)_H 
&=& - (\partial_t \varphi_{\varepsilon} (s), \varphi_{\varepsilon}(s)-m_0 )_{V^*} \\ \nonumber
&\leq& \| \partial_t \varphi_{\varepsilon} (s)\|_{V^*} (\| \varphi_{\varepsilon} (s)\|_{V^*} + |m_0| |\Omega|) \\ \label{stimaoperatoreN}
&=& c \| \partial_t \varphi_{\varepsilon}(s)\|_{V^*} (\| \varphi_{\varepsilon} (s)\|_H + 1). 
\end{eqnarray}
Recalling \eqref{eluiilveroriferimentoresatto}, we have that
\begin{equation} \label{stimamu}
(m(\mu_{\varepsilon}(s)), \varphi_{\varepsilon}(s)-m_0)_H = m(\mu_{\varepsilon}(s)) \Bigg( \int_{\Omega} \varphi_{\varepsilon}(s)-| \Omega | m_0\Bigg) = 0.
\end{equation}
Due to the Neumann homogeneous boundary conditions for $\varphi_{\varepsilon}$, we have that
\begin{equation} \label{precedent}
\int_{\Omega} \Delta \varphi_{\varepsilon}(s) =0.
\end{equation}
Thanks to \eqref{precedent}, we infer that
\begin{equation} \label{nustima}
\nu (\Delta \varphi_{\varepsilon}(s), \varphi_{\varepsilon}(s)-m_0)_H = - \nu \| \nabla \varphi_{\varepsilon}(s)\|^2_H - m_0 \int_{\Omega} \Delta \varphi_{\varepsilon}(s) =  - \nu \| \nabla \varphi_{\varepsilon}(s)\|^2_H \leq 0.
\end{equation}
As $\pi$ is a Lipschitz continuous function with Lipschitz constant $C_{\pi}$, we obtain that
\begin{eqnarray} \nonumber
- (\pi(\varphi_{\varepsilon}(s)), \varphi_{\varepsilon}(s)-m_0)_H 
&\leq& \int_{\Omega} |\pi(\varphi_{\varepsilon}(s))|| \varphi_{\varepsilon}(s)-m_0| \\ \nonumber
&\leq& \int_{\Omega} \bigg( |\pi(\varphi_{\varepsilon}(s)) - \pi(0)| + |\pi(0)| \bigg) \bigg(|\varphi_{\varepsilon}(s)|+ |m_0|\bigg) \\ \nonumber
&\leq& \int_{\Omega} \bigg( C_{\pi}|\varphi_{\varepsilon}(s)| + |\pi(0)| \bigg) \bigg(|\varphi_{\varepsilon}(s)|+ |m_0|\bigg) \\ \nonumber
&\leq& C_{\pi} \| \varphi_{\varepsilon}(s) \|^2_H 
+ \bigg(C_{\pi}|m_0| + |\pi(0)|\bigg)\| \varphi_{\varepsilon}(s) \|^2_H + c  \\ \label{stimaplemamma1}
&\leq& c (  \| \varphi_{\varepsilon}(s) \|^2_H + 1).
\end{eqnarray}
Moreover, we have that
\begin{eqnarray} \nonumber
\gamma (\chi_{\varepsilon}(s), \varphi_{\varepsilon}(s)-m_0)_H 
&\leq & \gamma \int_{\Omega} |\chi_{\varepsilon}(s)|| \varphi_{\varepsilon}(s)| + \gamma |m_0| \int_{\Omega}  |\chi_{\varepsilon}(s)| \\ \nonumber
&\leq & \gamma \|\chi_{\varepsilon}(s)\|^2_H + \gamma \| \varphi_{\varepsilon}(s)\|^2_H + \gamma |m_0| \|\chi_{\varepsilon}(s)\|^2_H + \gamma |m_0| |\Omega| \\ \label{stimagammateta}
&\leq & c (\|\chi_{\varepsilon}(s)\|^2_H +  \| \varphi_{\varepsilon}(s)\|^2_H + 1).
\end{eqnarray}
Consequently, rearranging the right-hand side of \eqref{astrodiriferimento} using \eqref{stimaoperatoreN}--\eqref{stimamu} and \eqref{nustima}--\eqref{stimagammateta}, we obtain that
\begin{equation} \label{astrodiriferimento22}
(\xi_{\varepsilon}(s), \varphi_{\varepsilon}(s)-m_0)_H
\leq c \bigg( \| \partial_t \varphi_{\varepsilon}(s)\|_{V^*} 
+  \| \varphi_{\varepsilon} (s)\|^2_H 
+  \|\chi_{\varepsilon}(s)\|^2_H + 1 \bigg) .
\end{equation}
Due to a useful inequality stated in \cite[Section 5]{GMS09}, it turns out that 
\begin{equation} \label{stimagilardi}
| \xi_{\varepsilon}(s)| \leq c[ \xi_{\varepsilon}(s) (\varphi_{\varepsilon}(s) -m_0) + 1].
\end{equation}
We integrate \eqref{stimagilardi} over $\Omega$ and,
due to \eqref{astrodiriferimento22}, we infer that
\begin{equation} \nonumber
\| \xi_{\varepsilon}(s) \|_{L^1(\Omega)} 
\leq c \bigg[ (\xi_{\varepsilon}(s), \varphi_{\varepsilon}(s) - m_0)_H + 1 \bigg]
\end{equation}
\begin{equation} \label{ultimalemma2222}
\leq c \bigg( \| \partial_t \varphi_{\varepsilon} (s)\|_{V^*} 
+  \| \varphi_{\varepsilon} (s)\|^2_H 
+  \|\chi_{\varepsilon}(s)\|^2_H +  1 \bigg).
\end{equation}
Due to \eqref{primastimaaprioiri1}--\eqref{primastimaaprioiri3}, from \eqref{ultimalemma2222}
we conclude that there exists a positive constant $c$, independent of $\varepsilon$, such that
\begin{equation} \label{sintesilemma2versionestretta}
\| \xi_{\varepsilon} \|_{L^2(0,T;L^1(\Omega))} \leq c.
\end{equation}

\paragraph{Third a priori estimate.}
As $\pi$ is a Lipschitz continuous function with Lipschitz constant
$C_{\pi}$, for every $s \in (0, T)$ we have that
\begin{eqnarray} \nonumber
|\pi(\varphi_{\varepsilon}(s))|^2 
&\leq& \big( |\pi(\varphi_{\varepsilon}(s)) - \pi(0)| + | \pi(0)| \big)^2  \\ \nonumber
&\leq& \big( C_{\pi} |\varphi_{\varepsilon}(s)| + | \pi(0)| \big)^2  \\ \label{utileulterioirestimasupigreco}
&\leq& c \big( |\varphi_{\varepsilon}(s)|^2 + 1 \big). 
\end{eqnarray}
Now, integrating \eqref{riferimentostima5} over 
$\Omega$, squaring the resultant and using \eqref{primastimaaprioiri1}--\eqref{primastimaaprioiri2} 
and \eqref{utileulterioirestimasupigreco}, we obtain that
\begin{equation} \nonumber
|m(\mu_{\varepsilon}(s))|^2 \leq  \frac{3}{|\Omega|^2} 
\bigg( \| \xi_{\varepsilon}(s) \|^2_{L^1(\Omega)} 
+ |\Omega| \| \pi(\varphi_{\varepsilon}(s))\|^2_H
+ \gamma \| \chi_{\varepsilon}(s) \|^2_H\bigg) 
\end{equation}
\begin{equation}  \label{finalelemma3}
\leq c
\bigg( \| \xi_{\varepsilon}(s) \|^2_{L^1(\Omega)} 
+ \| \varphi_{\varepsilon}(s) \|_H^2 
+ \| \chi_{\varepsilon}(s) \|^2_H + 1 \bigg).
\end{equation}
Consequently, integrating
\eqref{finalelemma3} over $(0,T)$ and recalling the previous a priori estimates \eqref{primastimaaprioiri1}--\eqref{primastimaaprioiri3} and \eqref{sintesilemma2versionestretta}, we conclude that there exists a positive constant 
$c$, independent of $\varepsilon$, such that
\begin{equation} \label{sintesilemma423333333333333333333333333333333}
\| m(\mu_{\varepsilon})\|_{L^2(0,T)}  \leq c.
\end{equation}

\paragraph{Fourth a priori estimate.}
We recall that the Poincar$\acute{e}$ inequality states that there exists a positive
constant $c_p$ such that
\begin{equation} \label{poincare}
\|z\|^2_V \leq c_p \| \nabla z \|^2_H \ \ \ \textrm{for all $z \in V$ with $m(z)=0$}.
\end{equation}
We integrate over $(0,T)$ the square of the norms in $V$
of each term of \eqref{riferimentostima4}.
Then, applying \eqref{sintesilemma423333333333333333333333333333333} and \eqref{poincare}, we obtain that
\begin{eqnarray} \nonumber
\int_0^T \| \mu_{\varepsilon}(s) \|^2_V \ ds &\leq& 2 \int_0^T \| m(\mu_{\varepsilon}(s)) \|^2_V \ ds 
+ 2 \int_0^T \|  \mathcal{N}\partial_t \varphi_{\varepsilon} (s) \|^2_V \ ds \\ \nonumber
&\leq& 2 \int_0^T | m(\mu_{\varepsilon}(s)) |^2 \ ds 
+ 2 c_p \int_0^T \| \nabla \mathcal{N}\partial_t \varphi_{\varepsilon} (s) \|^2_H  \ ds \\
&\leq& c + 2 c_p \int_0^T \| \partial_t \varphi_{\varepsilon} (s) \|^2_{V^*} \ ds .
\end{eqnarray}
Due to \eqref{primastimaaprioiri3}, we conclude that 
there exists a positive constant $c$, independent of $\varepsilon$, such that
\begin{equation} \label{sintesilemma42}
\| \mu_{\varepsilon} \|_{L^2(0,T;V)} \leq c.
\end{equation}

\paragraph{Fifth a priori estimate.}
We test \eqref{riferimentostima3} at time $s \in (0,T)$
by $\xi_{\varepsilon}(s) \in V$ and integrate the resultant 
over $\Omega$. We obtain that
\begin{equation} \label{primalemma5}
\| \xi_{\varepsilon}(s) \|^2_H 
= \big( \mu_{\varepsilon}(s) + \nu \Delta \varphi_{\varepsilon}(s) - \pi(\varphi_{\varepsilon}(s)) 
+ \gamma \chi_{\varepsilon}(s), \xi_{\varepsilon}(s) \big)_H. 
\end{equation}
Due to the monotonicity of $\beta_{\varepsilon}$, we have that
\begin{eqnarray} \nonumber
(\nu \Delta \varphi_{\varepsilon}(s), \xi_{\varepsilon}(s))_H 
&=& \nu \int_{\Omega} \Delta \varphi_{\varepsilon}(s) \xi_{\varepsilon}(s) \\ \nonumber
&=& - \nu \int_{\Omega} \nabla \varphi_{\varepsilon}(s) \cdot \nabla \xi_{\varepsilon}(s) \\ \label{utilesimapezzolemma5}
&=& - \nu \int_{\Omega} |\nabla \varphi_{\varepsilon}(s)|^2 \beta_{\varepsilon}'(\varphi_{\varepsilon}(s)) \leq 0.
\end{eqnarray}
Using \eqref{utilesimapezzolemma5} and the Young inequality, we can estimate \eqref{primalemma5} as follows
\begin{eqnarray} \nonumber
\| \xi_{\varepsilon}(s) \|^2_H 
& \leq & \big( \mu_{\varepsilon}(s) - \pi(\varphi_{\varepsilon}(s)) 
+ \gamma \chi_{\varepsilon}(s), \xi_{\varepsilon}(s)\big)_H \\ \nonumber
& \leq & \|\mu_{\varepsilon}(s) - \pi(\varphi_{\varepsilon}(s)) 
+ \gamma \chi_{\varepsilon}(s) \|_H \| \xi_{\varepsilon}(s) \|_H \\ \label{ultimosforzo}
& \leq & \frac{1}{2} \| \xi_{\varepsilon}(s) \|^2_H
+ 2 \big( \| \mu_{\varepsilon}(s)\|^2_H + \| \pi(\varphi_{\varepsilon}(s)) \|^2_H
+ \gamma^2 \| \chi_{\varepsilon}(s) \|^2_H \big).
\end{eqnarray}
Due to \eqref{utileulterioirestimasupigreco}, from \eqref{ultimosforzo} we infer that
\begin{equation} \label{passaggiochiave5}
\| \xi_{\varepsilon}(s) \|^2_H 
\leq  c \big( \| \mu_{\varepsilon}(s)\|^2_H + \| \varphi_{\varepsilon}(s) \|^2_H 
+ \| \chi_{\varepsilon}(s) \|^2_H + 1 \big).
\end{equation}
Then, integrating \eqref{passaggiochiave5}
over $(0,T)$ with respect to $s$ and using 
\eqref{primastimaaprioiri1}--\eqref{primastimaaprioiri3} and
\eqref{sintesilemma42}, we have that
\begin{equation} \label{sintesilemma51}
\| \xi_{\varepsilon} \|_{L^2(0,T;H)} \leq c,
\end{equation}
for some positive constant $c$, independent of $\varepsilon$.

\paragraph{Sixth a priori estimate.}
We integrate over $(0,T)$ the square of the norms in $H$
of each term of  \eqref{riferimentostima3}. 
Then, using \eqref{utileulterioirestimasupigreco}, \eqref{sintesilemma42} and \eqref{sintesilemma51}, we obtain that
\begin{equation} \nonumber
\nu^2 \int_0^T \| \Delta \varphi_{\varepsilon}(s) \|^2_H \ ds 
\end{equation} 
\begin{equation} \nonumber
\leq 4\int_0^T \| \mu_{\varepsilon}(s)  \|^2_H \ ds
+4 \int_0^T  \| \xi_{\varepsilon}(s) \|^2_H \ ds
+4 \int_0^T \|\pi(\varphi_{\varepsilon}(s)) \|^2_H \ ds
+4 \gamma^2 \int_0^T \| \chi_{\varepsilon}(s) \|^2_H \ ds 
\end{equation}
\begin{equation}
\leq c \Bigg(\int_0^T \| \varphi_{\varepsilon}(s) \|^2_H \ ds  
+ \int_0^T \| \chi_{\varepsilon}(s) \|^2_H \ ds + 1 \Bigg).
\end{equation}
Thanks to \eqref{primastimaaprioiri1}--\eqref{primastimaaprioiri3}, we conclude that there exists a positive constant $c$, independent of $\varepsilon$, such that
\begin{equation} \label{sintesilemma62}
\| \varphi_{\varepsilon} \|_{ L^2(0,T;W)} \leq c.
\end{equation}

\paragraph{Summary of the a priori estimates.}
Let us summarize the a priori estimates. From \eqref{primastimaaprioiri1}--\eqref{primastimaaprioiri2}, 
\eqref{sintesilemma42}, \eqref{sintesilemma51} and \eqref{sintesilemma62}
we conclude that there exists a constant $c>0$, independent of $\varepsilon$,
such that
\begin{eqnarray} \label{stimaapriorisintesi1}
\| \chi_{\varepsilon} \|_{H^1(0,T; V^*) \cap L^{\infty}(0,T; H) \cap L^2(0,T; V)} & \leq  & c  , \\ \label{stimaapriorisintesi2}
\| \varphi_{\varepsilon} \|_{H^1(0,T; V^*) \cap L^{\infty}(0,T; V) \cap L^2(0,T;W)}  &\leq &  c, \\ \label{stimaapriorisintesi3}
\| \zeta_{\varepsilon} \|_{ L^{\infty}(0,T; H)} & \leq & c, \\ \label{stimaapriorisintesi4}
\| \xi_{\varepsilon} \|_{L^2(0,T;H)} &\leq & c , \\ \label{stimaapriorisintesi5}
\| \mu_{\varepsilon} \|_{L^2(0,T;V)} &\leq & c .
\end{eqnarray}

\section{Existence - Passage to the limit as $\varepsilon \searrow 0$}
\setcounter{equation}{0}
Based on available results (cf., e.g., \cite{CF15a}), it turns out that
there exists a solution $(\chi_{\varepsilon},\varphi_{\varepsilon},\mu_{\varepsilon})$ of ($P_{\varepsilon}$)
satisfying the regularity requirements \eqref{regolaritaesistenza-epsilon1}--\eqref{regolaritaesistenza-epsilon3}
and solving \eqref{iniziale1-epsilon}-\eqref{iniziale7-epsilon}.
In this section we pass to the limit as $\varepsilon \searrow 0$ and prove that
the limit of subsequences of solutions $(\chi_{\varepsilon},\varphi_{\varepsilon},\mu_{\varepsilon})$
for ($P_{\varepsilon}$) (see \eqref{iniziale1-epsilon}--\eqref{iniziale7-epsilon}) yields a solution $(\chi,\varphi,\mu)$
of ($P$) (see \eqref{iniziale1}--\eqref{iniziale7}).

Thanks to the uniform estimates \eqref{stimaapriorisintesi1}--\eqref{stimaapriorisintesi5},
there exists a subsequence $\{ \varepsilon _k \}_{k \in \mathbb{N}}$ with 
$\varepsilon _k \searrow  0$ as $k \to +\infty $
and some limit functions $\chi \in H^1(0,T;V^*) \cap L^\infty (0,T;H) \cap L^2(0,T; V)$,
$\varphi \in H^1(0,T;V^*) \cap L^\infty (0,T;H) \cap L^2(0,T;W)$,
$\mu \in L^2(0,T;V)$,
$\xi \in L^2(0,T;H)$
and $\zeta \in L^\infty (0,T;H)$ such that 
\begin{eqnarray} \label{limitedebole1}
\chi_{\varepsilon_k}  \rightharpoonup^* \chi  & \textrm{in} & H^1(0,T; V^*) \cap L^{\infty}(0,T; H), \\ \label{limitedebole2}
\chi_{\varepsilon_k}  \rightharpoonup \chi    & \textrm{in} & L^2(0,T; V), \\ \label{limitedebole3}
\varphi_{\varepsilon_k}    \rightharpoonup^* \varphi    & \textrm{in} & H^1(0,T; V^*) \cap L^{\infty}(0,T; V), \\ \label{limitedebole4}
\varphi_{\varepsilon_k}    \rightharpoonup \varphi      & \textrm{in} & L^2(0,T;W), \\ \label{limitedebole5}
\mu_{\varepsilon_k}        \rightharpoonup   \mu        & \textrm{in} & L^2(0,T;V), \\ \label{limitedebole6}
\xi_{\varepsilon_k}        \rightharpoonup   \xi        & \textrm{in} & L^2(0,T;H), \\ \label{limitedebole7}
\zeta_{\varepsilon_k}       \rightharpoonup^*  \zeta      & \textrm{in} & L^{\infty} (0,T;H),
\end{eqnarray} 
as $ k \to + \infty$. 
From \eqref{limitedebole1}--\eqref{limitedebole4} and the well-known Ascoli--Arzel\'a theorem
(see, e.g., \cite[Sect. 8, Cor. 4]{Sim87}), we infer that 
\begin{eqnarray} \label{limiteforte1} 
\chi_{\varepsilon_k}  \longrightarrow \chi  & \textrm{in} & C^0([0,T];V^*) \cap L^{2}(0,T; H), \\ \label{limiteforte2}
\varphi_{\varepsilon_k}    \longrightarrow \varphi    & \textrm{in} & C^0([0,T];H) \cap L^{2}(0,T; V),
\end{eqnarray} 
as  $ k \to + \infty$. 
As $\pi$ is a Lipschitz continuous function, 
for a.e. $s \in [0,T]$ we have that
\begin{equation} \nonumber
|\pi(\varphi_{\varepsilon_k}(s)) - \pi(\varphi(s))|
\leq  C_{\pi} |\varphi_{\varepsilon_k}(s) - \varphi(s)|.
\end{equation}
Thanks to \eqref{limiteforte2}, we conclude that
\begin{equation} \label{convergenzaepsilondipigreco}
\pi(\varphi_{\varepsilon_k}(s)) \longrightarrow \pi(\varphi(s)) \quad \textrm{in $L^2(0,T; H)$,}
\end{equation}
as  $ k \to + \infty$. 

\paragraph{Passage to the limit on $\xi_{\varepsilon}$.}
In this paragraph we check that $\xi \in \beta(\varphi)$ a.e. in $Q$.
To this aim, we recall that
\begin{equation} \label{deboleepsilon1111}
\varphi_{\varepsilon_k} \rightarrow \varphi \quad \quad \textrm{in} \  L^2(0,T;H) \equiv L^2(Q),
\end{equation}
\begin{equation} \label{deboleepsilon2222}
\xi_{\varepsilon_k} \rightharpoonup \xi \quad \quad \textrm{in} \  L^2(0,T;H),
\end{equation}
as $ k \to + \infty$. Now, we introduce the operator $\mathcal{B}_{\varepsilon}$ induced by $\beta_{\varepsilon}$ on $L^2(Q)$ in the following way
\begin{equation} \nonumber
\mathcal{B}_{\varepsilon}: L^2(Q) \longrightarrow L^2(Q)
\end{equation}
\begin{equation} 
\xi_{\varepsilon}  \in \mathcal{B}_{\varepsilon}(\varphi_{\varepsilon}) \Longleftrightarrow \xi_{\varepsilon}(x,t) \in \beta_{\varepsilon}(\varphi_{\varepsilon}(x,t)) \quad \textrm{for a.e. $(x,t) \in Q$.}
\end{equation}
Due to \eqref{deboleepsilon1111}--\eqref{deboleepsilon2222}, as $ k \to + \infty$, we have that 
\begin{equation} \label{fine11}
\left\{ \begin{array}{ll}
\mathcal{B}_{\varepsilon_k}(\varphi_{\varepsilon_k}) \rightharpoonup \xi & \textrm{in $L^2(Q)$,} \\
\varphi_{\varepsilon_k} \rightarrow \varphi & \textrm{in $L^2(Q)$,}
\end{array}
\right.
\end{equation}
\begin{equation} \label{fine12}
\limsup_{k \rightarrow +\infty} \int_Q \xi_{\varepsilon_k} \varphi_{\varepsilon_k} = \int_Q \xi \varphi.
\end{equation}
Thanks to \eqref{fine11}--\eqref{fine12} and to the general result \cite[Proposition~2.2, p.~38]{Bar10}, we conclude that
\begin{equation} \label{proprioquestaqui}
\xi  \in \mathcal{B}(\varphi) \quad \textrm{in $L^2(Q)$},
\end{equation}
with analogous definition for $\mathcal{B}$ (see \eqref{betagrande1}--\eqref{betagrande2}). This
is equivalent to saying that
\begin{equation} \label{ultimaconvergenzadibeta} 
\xi  \in \beta(\varphi) \quad \textrm{a.e. in $Q$}. 
\end{equation}

\paragraph{Passage to the limit on $\zeta_{\varepsilon}$.}
In this paragraph we check that $\zeta(t) \in A(a \chi(t) + b \varphi(t) - \kappa^*)$ for a.e. $t \in [0,T]$.
Let us recall that
\begin{eqnarray} \label{siquestaserve1}
\chi_{\varepsilon_k} \rightarrow \chi  & \textrm{in} & L^2 (0,T;H), \\ \label{siquestaserve2}
\varphi_{\varepsilon_k}   \rightarrow \varphi    & \textrm{in} & L^2 (0,T;H), \\ \label{siquestaserve2222222}
\zeta_{\varepsilon_k}   \rightharpoonup \zeta    & \textrm{in} & L^2 (0,T;H),
\end{eqnarray} 
as $ k \to + \infty$. Setting
\begin{equation} \nonumber
\kappa_{\varepsilon_k} := a \chi_{\varepsilon_k} + b\varphi_{\varepsilon_k} - \kappa^*, \quad \quad
\kappa := a \chi + b\varphi - \kappa^*,
\end{equation}
thanks to \eqref{siquestaserve1}--\eqref{siquestaserve2}, we have that
\begin{equation}
\kappa_{\varepsilon_k} \longrightarrow \kappa \quad \textrm{in} \ L^{2} (0,T;H),
\end{equation}
as $ k \to + \infty$. Now, we introduce the operator $\mathcal{A}_{\varepsilon}$ induced by $A_{\varepsilon}$ on $L^2(0,T;H)$ in the following way
\begin{equation} \nonumber
\mathcal{A}_{\varepsilon}: L^2(0,T;H) \longrightarrow L^2(0,T;H)
\end{equation}
\begin{equation} 
\zeta_{\varepsilon}  \in \mathcal{A}_{\varepsilon}(\kappa_{\varepsilon}) \Longleftrightarrow \zeta_{\varepsilon}(t) \in A_{\varepsilon}(\kappa_{\varepsilon}(t)) \quad 
\textrm{for a.e. $t \in [0,T]$.}
\end{equation}
Due to \eqref{siquestaserve1}--\eqref{siquestaserve2222222}, we have that
\begin{equation} \label{fine21}
\left\{ \begin{array}{ll}
\mathcal{A}_{\varepsilon_k}(\kappa_{\varepsilon_k}) \rightharpoonup \zeta & \textrm{in $L^2(0,T;H)$,} \\
\kappa_{\varepsilon_k} \rightarrow \kappa & \textrm{in $L^2(0,T;H)$,}
\end{array}
\right.
\end{equation}
\begin{equation} \label{fine22}
\limsup_{k \rightarrow +\infty} \int_Q \zeta_{\varepsilon_k} \kappa_{\varepsilon_k} = \int_Q \zeta \kappa.
\end{equation}
Thanks to \eqref{fine21}--\eqref{fine22} and the convergence result \cite[Proposition~2.2, p.~38]{Bar10}, we conclude that
\begin{equation} 
\zeta \in \mathcal{A}(\kappa) \quad \textrm{in $L^2(0,T;H)$}, 
\end{equation}
with obvious definition for $\mathcal{A}$ (see \eqref{Agrande1}--\eqref{Agrande2}). This is equivalent to saying that
\begin{equation} \label{ultimaconvergenzadiA}
\zeta(t) \in A( a \chi(t) + b\varphi(t) - \kappa^* ) \quad \textrm{for a.e. $t \in [0,T]$.} 
\end{equation}
\paragraph{Conclusion of the proof}
Using \eqref{limitedebole1}--\eqref{convergenzaepsilondipigreco},   \eqref{ultimaconvergenzadibeta} and \eqref{ultimaconvergenzadiA},
we can pass to the limit as $\varepsilon \searrow 0$ in \eqref{iniziale1-epsilon}--\eqref{iniziale7-epsilon}
obtaining \eqref{iniziale1}--\eqref{iniziale7} for the limiting functions $\chi$, $\varphi$ and $\mu$.

\setcounter{equation}{0}
\section{Regularity}
This section is devoted to the proof of Theorem \ref{Regolar}. 
In order to obtain additional regularity for the solutions, we need further a priori estimates obtained from the approximating problem 
$(P_{\varepsilon})$ (see \eqref{iniziale1-epsilon}--\eqref{iniziale7-epsilon})
in which we take $\chi_{0\varepsilon} = \chi_0$ and $\varphi_{0\varepsilon}= \varphi_0$.

\paragraph{Seventh a priori estimate.}
We test \eqref{iniziale1-epsilon} by $\partial_t \chi_{\varepsilon}$ and integrate over $Q_t$, $t \in (0,T]$. We have that
\begin{equation} \label{reg1}
\int_{Q_t} |\partial_t\chi_{\varepsilon}|^2 + \ell \int_{Q_t} \partial_t \varphi_{\varepsilon} \partial_t \chi_{\varepsilon} 
+ \frac{1}{2} \int_{\Omega} |\nabla \chi_{\varepsilon}(t)| 
+ \int_{Q_t} \zeta_{\varepsilon} \partial_t \chi_{\varepsilon} 
= \int_{Q_t} f_{\varepsilon} \partial_t \chi_{\varepsilon} + \frac{1}{2} \int_{\Omega} |\nabla \chi_{0}| .
\end{equation}
We now proceed with a formal estimate since we have to differentiate \eqref{iniziale2-epsilon} and \eqref{iniziale3-epsilon}
with respect to time. For a rigorous approach, one can argue, e.g., as in \cite[Subsection 4.4]{CF15bbbbbb}.
By time differentiation of \eqref{iniziale2-epsilon} and \eqref{iniziale3-epsilon} we have
\begin{equation} \label{reg2}
\partial_{tt} \varphi_{\varepsilon} - \Delta \partial_t \mu_{\varepsilon} = 0  ,
\end{equation}
\begin{equation} \label{reg3}
\partial_t \mu_{\varepsilon}= - \nu \Delta \partial_t \varphi_{\varepsilon} + \beta_{\varepsilon}'(\varphi_{\varepsilon})\partial_t \varphi_{\varepsilon} + \pi'(\varphi_{\varepsilon})\partial_t \varphi_{\varepsilon} - \gamma \partial_t \chi_{\varepsilon}.
\end{equation}
According to \eqref{eluiilveroriferimentoresatto}, $m(\partial_t \varphi_{\varepsilon})=0$. Consequently,
$\partial_t \varphi_{\varepsilon} \in D(\mathcal{N})$ and we can test \eqref{reg2} by $\frac{\ell}{\gamma} \mathcal{N} (\partial_t \varphi_{\varepsilon})$.
Integrating the resultant over $Q_t$, we obtain that
\begin{equation} \label{Regg2}
- \frac{\ell}{\gamma} \int_{Q_t} \partial_t \mu_{\varepsilon} \partial_t \varphi_{\varepsilon} 
= \frac{\ell}{2\gamma}  \| \partial_t \varphi_{\varepsilon} (t) \|^2_{V^*}  
- \frac{\ell}{2\gamma}  \| \partial_t \varphi_{\varepsilon}(0) \|^2_{V^*}   .
\end{equation}
We test \eqref{reg3} by $\frac{\ell}{\gamma} \partial_t \varphi_{\varepsilon}$ and integrate over $Q_t$. We have that
\begin{equation} \nonumber
\frac{\ell}{\gamma} \int_{Q_t} \partial_t \mu_{\varepsilon} \partial_t \varphi_{\varepsilon}
\end{equation} 
\begin{equation} \label{Regg3}
= \frac{\nu \ell}{\gamma}  \int_{Q_t} |\nabla \partial_t \varphi_{\varepsilon}|^2 
+ \frac{\ell}{\gamma} \int_{Q_t} \beta_{\varepsilon}'(\varphi_{\varepsilon}) |\partial_t \varphi_{\varepsilon}|^2 
+ \frac{\ell}{\gamma} \int_{Q_t} \pi'(\varphi_{\varepsilon}) |\partial_t \varphi_{\varepsilon}|^2 
- \ell \int_{Q_t} \partial_t \varphi_{\varepsilon} \partial_t \chi_{\varepsilon}.
\end{equation}
By combining \eqref{reg1},  \eqref{Regg2} and \eqref{Regg3}, we infer that
\begin{equation} \nonumber
\int_{Q_t} |\partial_t\chi_{\varepsilon}|^2 
+ \frac{1}{2} \int_{\Omega} |\nabla \chi_{\varepsilon}(t)| 
+ \frac{\nu \ell}{\gamma}  \int_{Q_t} |\nabla \partial_t \varphi_{\varepsilon}|^2 
+ \frac{\ell}{2\gamma}  \| \partial_t \varphi_{\varepsilon} (t) \|^2_{V^*}
= \frac{1}{2} \int_{\Omega} |\nabla \chi_{0}| 
\end{equation}
\begin{equation} \label{stimaregol0}
+ \int_{Q_t} f_{\varepsilon} \partial_t \chi_{\varepsilon} 
+ \frac{\ell}{2\gamma}  \| \partial_t \varphi_{\varepsilon}(0) \|^2_{V^*}
- \int_{Q_t} \zeta_{\varepsilon} \partial_t \chi_{\varepsilon} 
- \frac{\ell}{\gamma} \int_{Q_t} \beta_{\varepsilon}'(\varphi_{\varepsilon}) |\partial_t \varphi_{\varepsilon}|^2 
- \frac{\ell}{\gamma} \int_{Q_t} \pi'(\varphi_{\varepsilon}) |\partial_t \varphi_{\varepsilon}|^2 .
\end{equation}
By applying inequality \eqref{dis2} to the second term on the right-hand side of \eqref{stimaregol0}, we infer that
\begin{equation}  \label{11}
\int_{Q_t} f_{\varepsilon} \partial_t \chi_{\varepsilon} \leq \| f_{\varepsilon} \|^2_{L^2(0,T;H)} + \frac{1}{4} \int_{Q_t} |\partial_t\chi_{\varepsilon}|^2 .
\end{equation}
Moreover, as $\beta_{\varepsilon}$ is a maximal monotone operator, we have that $\beta_{\varepsilon}'>0$
and consequently
\begin{equation} \label{22}
- \frac{\ell}{\gamma} \int_{Q_t} \beta_{\varepsilon}'(\varphi_{\varepsilon}) |\partial_t \varphi_{\varepsilon}|^2 \leq 0.
\end{equation}
Due to \eqref{primastimaaprioiri5}, we have that
\begin{equation} \label{44}
- \int_{Q_t} \zeta_{\varepsilon} \partial_t \chi_{\varepsilon} 
\leq  \int_{Q_t} |\zeta_{\varepsilon}|^2 +  \frac{1}{4} \int_{Q_t} |\partial_t \chi_{\varepsilon}|^2 
\leq c +  \frac{1}{4} \int_{Q_t} |\partial_t \chi_{\varepsilon}|^2.
\end{equation}
As $\pi$ is a Lipschitz continuous function with Lipschitz constant $C_{\pi}$,
we have that
\begin{equation} \label{55}
- \frac{\ell}{\gamma} \int_{Q_t} \pi'(\varphi_{\varepsilon}) |\partial_t \varphi_{\varepsilon}|^2 
\leq \frac{\ell}{\gamma} \int_{Q_t} |\pi'(\varphi_{\varepsilon})| |\partial_t \varphi_{\varepsilon}|^2
\leq \frac{C_{\pi}\ell}{\gamma} \int_{Q_t} |\partial_t \varphi_{\varepsilon}|^2.
\end{equation}
Adding $\frac{\nu \ell}{\gamma} \int_{Q_t} |\partial_t \varphi_{\varepsilon}|^2 $ 
to both side of \eqref{stimaregol0} and rearranging the right-hand side of \eqref{stimaregol0}
using \eqref{11}--\eqref{55},  we obtain that
\begin{equation} \nonumber
\frac{1}{2} \int_{Q_t} |\partial_t \chi_{\varepsilon}|^2  
+ \frac{1}{2} \int_{\Omega} |\nabla \chi_{\varepsilon}(t)| 
+ \frac{\nu \ell}{\gamma}  \int_0^t \| \partial_t \varphi_{\varepsilon}(s)\|_V^2 \ ds
+ \frac{\ell}{2\gamma}  \| \partial_t \varphi_{\varepsilon} (t) \|^2_{V^*}
\end{equation}
\begin{equation}  \label{utilissima11112}
\leq \frac{1}{2} \|\chi_{0}\|^2_V + \frac{\ell}{2\gamma}  \| \partial_t \varphi_{\varepsilon}(0) \|^2_{V^*} + \| f_{\varepsilon} \|^2_{L^2(0,T;H)} 
+ \Bigg( \frac{C_{\pi}\ell}{\gamma} + \frac{\nu \ell}{\gamma}  \Bigg) \int_0^t \| \partial_t \varphi_{\varepsilon} (s)\|_{V^*}^2 \ ds + c.
\end{equation}
Thanks to the compactness of the embedding $V \subset H \subset V^*$, the inequality stated by \cite[Lemma 8, p. 84]{Sim87}
ensures that, choosing
\begin{equation} \nonumber
\delta  = \Bigg(\frac{\nu \ell}{4\gamma} \Bigg( \frac{C_{\pi}\ell}{\gamma} + \frac{\nu \ell}{\gamma}  \Bigg)^{-1} \Bigg)^{\frac{1}{2}},
\end{equation}
we can estimate the fourth term on the right-hand side of \eqref{utilissima11112} as follows
\begin{equation} \label{stimasiusa}
\Bigg( \frac{C_{\pi}\ell}{\gamma} + \frac{\nu \ell}{\gamma}  \Bigg) \int_0^t \| \partial_t \varphi_{\varepsilon} (s)\|_H^2 \ ds
\leq \frac{\nu \ell}{2\gamma}  \int_0^t \| \partial_t \varphi_{\varepsilon} (s)\|_V^2 \ ds + c \int_0^t \| \partial_t \varphi_{\varepsilon} (s)\|^2_{V^*} \ ds.
\end{equation}
Due to \eqref{stimasiusa}, from \eqref{utilissima11112} we have that
\begin{equation} \nonumber
\frac{1}{2} \int_{Q_t} |\partial_t \chi_{\varepsilon}|^2  
+ \frac{1}{2} \int_{\Omega} |\nabla \chi_{\varepsilon}(t)| 
+ \frac{\nu \ell}{2\gamma}  \int_0^t \| \partial_t \varphi_{\varepsilon}(s)\|_V^2 \ ds
+ \frac{\ell}{2\gamma}  \| \partial_t \varphi_{\varepsilon} (t) \|^2_{V^*}
\end{equation}
\begin{equation} \nonumber
\leq \frac{1}{2} \|\chi_{0}\|^2_V + \frac{\ell}{2\gamma}  \| \partial_t \varphi_{\varepsilon}(0) \|^2_{V^*} 
+ \| f_{\varepsilon} \|^2_{L^2(0,T;H)} 
+ c  \int_0^t \frac{\ell}{2\gamma} \| \partial_t \varphi_{\varepsilon} (s)\|^2_{V^*} \ ds
\end{equation}
\begin{equation} \label{utilissima111122222}
\leq \frac{1}{2} \|\chi_{0}\|^2_V + \frac{\ell}{2\gamma}  \| \partial_t \varphi_{\varepsilon}(0) \|^2_{V^*} 
+ \| f_{\varepsilon} \|^2_{L^2(0,T;H)} 
+ c  \| \varphi_{\varepsilon} \|^2_{H^{1}(0,T; V^*)} +c.
\end{equation}
Since
$( -\nu \Delta \varphi_{0} + \beta_{\varepsilon}(\varphi_{0}) 
+ \pi(\varphi_{0}) - \gamma \chi_{0}) $
is bounded in $V$ uniformly with respect to $\varepsilon$ according to \eqref{tecnicadistima1},
we deduce, by comparison in \eqref{iniziale2-epsilon}--\eqref{iniziale3-epsilon}, 
that the second term on the right-hand side of \eqref{utilissima111122222} is estimated by a positive constant. 
Hence, due to \eqref{0reg}, \eqref{convergenzadatiinizialiepsilon1}--\eqref{regolaritaesistenza-epsilon3}
and \eqref{stimaapriorisintesi2},
the right-hand side of \eqref{utilissima111122222} is bounded and we conclude that there exists a positive constant $c$, 
independent of $\varepsilon$, such that
\begin{equation} \label{stimaregolulteriore1}
\| \chi_{\varepsilon} \|_{H^1(0,T;H) \cap L^\infty (0,T;V)} +   \| \varphi_{\varepsilon} \|_{W^{1, \infty}(0,T;V^*) \cap H^1(0,T;V)} \leq c .
\end{equation}

\paragraph{Eighth a priori estimate.}
From \eqref{iniziale1-epsilon}, we have that
\begin{equation} 
\Delta \chi_{\varepsilon} = \partial_t(\chi_{\varepsilon} + \ell \varphi_{\varepsilon}) + \zeta_{\varepsilon} - f_{\varepsilon} =: h_{\varepsilon}.
\end{equation}
We observe that \eqref{stimaregolulteriore1} ensures that $h_{\varepsilon}$ is bounded in $ {L^2 (0,T;H)}$
uniformly with respect to $\varepsilon$. Then we infer that
there exists a constant $c > 0$, independent of $\varepsilon$, such that
\begin{equation} \label{stimaregolulteriore2}
\| \chi_{\varepsilon} \|_{H^1(0,T;H) \cap L^\infty (0,T;V) \cap L^2 (0,T;W)}  \leq c .
\end{equation}

\paragraph{Ninth a priori estimate.}
Due to \eqref{stimaregolulteriore1}--\eqref{stimaregolulteriore2}, from \eqref{ultimalemma2222} we deduce that
\begin{equation} \label{aaa1}
\| \xi_{\varepsilon} \|_{ L^\infty (0,T;L^1(\Omega))}  \leq c .
\end{equation}
Now, using \eqref{finalelemma3}, we infer that $ \| m(\mu_{\varepsilon}) \|_{ L^\infty (0,T)}  \leq c $.
By comparison in \eqref{regolaritaesistenza-epsilon2} and \eqref{riferimentostima4}, it follows that
\begin{equation} \label{aaa3}
\| \mu_{\varepsilon} \|_{ L^\infty (0,T;V)}  \leq c .
\end{equation}
Moreover, from \eqref{passaggiochiave5}, we obtain that $\| \xi_{\varepsilon} \|_{ L^\infty (0,T;H)}  \leq c$.
Then, by comparison in \eqref{iniziale3-epsilon}, we conclude that
\begin{equation} \label{aaa5}
\| \Delta \varphi_{\varepsilon} \|_{ L^\infty (0,T;H) \cap L^2 (0,T;W)}  \leq c.
\end{equation}

\paragraph{Conclusion of the proof.}
As  \eqref{stimaregolulteriore1}, \eqref{stimaregolulteriore2} and \eqref{aaa1}--\eqref{aaa5} follow uniformly
with respect to $\varepsilon$, the same estimates hold true for the limiting functions $\chi$, $\varphi$ and $\mu$.
Hence, \eqref{regolaritaesistenza1reg}--\eqref{regolaritaesistenza3reg} are fulfilled and 
\begin{equation} 
\| \chi \|_{ H^1(0,T;H) \cap L^\infty (0,T;V) \cap  L^2(0,T;W) } \leq c, 
\end{equation}
\begin{equation} 
\| \varphi    \|_{ W^{1, \infty}(0,T;V^*) \cap H^1(0,T;V) \cap L^\infty (0,T;W)} \leq c, 
\end{equation}
\begin{equation} 
\| \mu        \|_{L^{\infty}(0,T;V)} \leq c.
\end{equation}

\setcounter{equation}{0}
\section{Uniqueness and continuous dependence}
This section is devoted to the proof of Theorem \ref{Teorema-unicita}.

Assume $a \ell = b$. If $f_i$, $\kappa^{*}_i$, $\kappa_{0_i}$, $\varphi_{0_i}$, $i=1,2$, 
are given as in \eqref{f}--\eqref{0} and $(\kappa_i, \varphi_i)$, $i=1,2$,  are the corresponding solutions of problem ($\widetilde{P}$) (see \eqref{iniziale1-uni1}--\eqref{iniziale7-uni1}), then we can write problem ($\widetilde{P}$) for both $(\kappa_i, \varphi_i)$, $i=1,2$
and take the difference between the respective equations.
Setting $\kappa:= \kappa_1 - \kappa_2$, $\varphi:=\varphi_1 - \varphi_2$,
$\mu := \mu_1 - \mu_2$, $f:= f_1 - f_2$, $\kappa^* := \kappa^*_1 - \kappa^*_2$,
$\kappa_0:= \kappa_{0_1} - \kappa_{0_2}$, $\varphi_0:= \varphi_{0_1} - \varphi_{0_2}$,
we obtain that
\begin{equation} \label{111}
\partial_t\kappa - \Delta \kappa + b \Delta \varphi - \Delta \kappa^* + a (\zeta_1 - \zeta_2) = a f, 
\end{equation}
\begin{equation} \label{222}
\partial_t \varphi - \Delta \mu= 0,  
\end{equation}
\begin{equation} \label{333}
\mu= - \nu \Delta \varphi + \xi_1 - \xi_2 + \pi(\varphi_1) - \pi(\varphi_2) - \frac{\gamma}{a} (\kappa - b \varphi + \kappa^*). 
\end{equation}
We observe that, due to \eqref{masseuguali}, $m(\varphi_0)=0$. 
Consequently, thanks to \eqref{integraleconservato11}, $m(\varphi)= 0$
and $\varphi \in D(\mathcal{N})$ a.e. in $(0,T)$ (see \eqref{operatoreNNNN1}).
Now, we test \eqref{111} by $\kappa$. Integrating over $Q_t$, $t \in (0,T]$, we have that
\begin{equation} \nonumber
\frac{1}{2} \int_{\Omega} |\kappa(t)|^2 + \int_{Q_t} |\nabla \kappa|^2 
- b \int_{Q_t} \nabla \varphi \cdot \nabla \kappa  + a \int_{Q_t} (\zeta_1 - \zeta_2) (\kappa_1 - \kappa_2)
\end{equation}
\begin{equation} \label{1primauni}
=  \frac{1}{2} \int_{\Omega} |\kappa_0|^2 + \int_{Q_t} (a f + \Delta \kappa^*)\kappa.
\end{equation}
We test \eqref{222} by $\frac{b^2}{\nu} \mathcal{N} \varphi$. Integrating over $(0,t)$, we obtain that
\begin{equation} \nonumber
\frac{b^2}{\nu} \int_0^t \langle \partial_t \varphi(s), \mathcal{N} \varphi(s) \rangle_{V^*,V} \ ds + \frac{b^2}{\nu}\int_{Q_t} \nabla \mu \cdot \nabla \mathcal{N} \varphi 
= 0, 
\end{equation} 
\begin{equation}  \label{2primauni}
\frac{b^2}{2 \nu} \| \varphi(t) \|^2_{V^*} + \frac{b^2}{\nu} \int_{Q_t} \mu \varphi = \frac{b^2}{2 \nu} \| \varphi_0 \|^2_{V^*}.
\end{equation} 
Testing \eqref{333} by $- \frac{b^2}{\nu}\varphi$ and integrating over $Q_t$, we have that
\begin{equation} \nonumber
-\frac{b^2}{\nu} \int_{Q_t} \mu \varphi = - b^2 \int_{Q_t} |\nabla \varphi|^2 
- \frac{b^2}{\nu} \int_{Q_t} (\xi_1 - \xi_2)(\varphi_1 - \varphi_2) 
\end{equation}
\begin{equation} \label{3primauni}
- \frac{b^2}{\nu} \int_{Q_t} [\pi(\varphi_1) - \pi(\varphi_2)](\varphi_1 - \varphi_2)
+ \frac{\gamma b^2}{a\nu} \int_{Q_t} (\kappa - b \varphi + \kappa^*) \varphi. 
\end{equation}
Then, we combine \eqref{1primauni}--\eqref{3primauni} and  infer that
\begin{equation} \nonumber
\frac{1}{2} \|\kappa(t)\|_H^2 + \int_{Q_t} ( |\nabla \kappa|^2 - b  \nabla \varphi \cdot \nabla \kappa +  b^2 |\nabla \varphi|^2 )
+ \frac{b^2}{2 \nu} \| \varphi(t) \|^2_{V^*}
\end{equation}
\begin{equation} \nonumber
+ a \int_{Q_t} (\zeta_1 - \zeta_2) (\kappa_1 - \kappa_2) + \frac{b^2}{\nu} \int_{Q_t} (\xi_1 - \xi_2)(\varphi_1 - \varphi_2) 
\end{equation}
\begin{equation} \nonumber
= - \frac{b^2}{\nu} \int_{Q_t} [\pi(\varphi_1) - \pi(\varphi_2)](\varphi_1 - \varphi_2)
+ \frac{\gamma b^2}{a\nu} \int_{Q_t} (\kappa - b \varphi + \kappa^*) \varphi
\end{equation}
\begin{equation} \label{riferimentounicita}
+ \frac{b^2}{2 \nu} \| \varphi_0 \|^2_{V^*} +  \frac{1}{2} \|\kappa_0\|_H^2 + \int_{Q_t} (a f + \Delta \kappa^*)\kappa.
\end{equation}
Since  $A$ and $\beta$ are maximal monotone, we have that
\begin{equation} \label{questa0}
a \int_{Q_t} (\zeta_1 - \zeta_2) (\kappa_1 - \kappa_2) \geq 0,  
\end{equation}
\begin{equation} \label{questa1}
\frac{b^2}{\nu} \int_{Q_t} (\xi_1 - \xi_2)(\varphi_1 - \varphi_2) \geq 0.
\end{equation}
Moreover, thanks to the Lipschitz continuity of $\pi$, we infer that 
\begin{eqnarray}  \nonumber
- \frac{b^2}{\nu} \int_{Q_t} [\pi(\varphi_1) - \pi(\varphi_2)](\varphi_1 - \varphi_2)
&\leq & \frac{b^2}{\nu} \int_{Q_t} |\pi(\varphi_1) - \pi(\varphi_2)||\varphi_1 - \varphi_2| \\ \label{questa2}
& \leq & \frac{C_{\pi} b^2}{\nu} \int_{Q_t} |\varphi|^2 .
\end{eqnarray} 
We also notice that the integral involving the gradients is estimated from below in this
way:
\begin{equation} \label{questa3}
\int_{Q_t} ( |\nabla \kappa|^2 - b  \nabla \varphi \cdot \nabla \kappa +  b^2 |\nabla \varphi|^2 ) \geq \frac{1}{2}
\int_{Q_t} ( |\nabla \kappa|^2 +  b^2 |\nabla \varphi|^2 ) .
\end{equation} 
Recalling that
\begin{equation} \label{stimissima333333}
- \frac{\gamma b^3}{a\nu} \int_{Q_t} |\varphi|^2 \leq 0,
\end{equation}
applying inequality \eqref{dis2} to the second and fifth term on the right-hand side of \eqref{riferimentounicita},
using \eqref{questa0}--\eqref{questa3} and adding to both sides $ b^2 \int_0^t \| \varphi(s)\|^2_H \ ds$, we infer that
\begin{equation} \nonumber
\frac{1}{2} \|\kappa(t)\|_H^2 + \int_{Q_t}  |\nabla \kappa|^2 +  b^2 \int_0^t \| \varphi(s)\|_V^2 \ ds 
+ \frac{b^2}{2 \nu} \| \varphi(t) \|^2_{V^*}
\end{equation}
\begin{equation} \label{precedentissima3}
\leq 
(K +b^2) \int_0^t \|\varphi(s)\|_H^2 \ ds  
+ \frac{1}{2} \int_{Q_t} |\kappa|^2
+ \frac{b^2}{2 \nu} \| \varphi_0 \|^2_{V^*}
+  \frac{1}{2} \|\kappa_0\|_H^2
+ 2 a^2 \|f\|_{L^2(0,T; H)}^2 + 3T \| \kappa^*\|_W^2,
\end{equation}
where
\begin{equation} \nonumber
K = \Bigg[ \frac{C_{\pi} b^2}{\nu} + 2 \Big( \frac{\gamma b^2}{a\nu} \Big)^2   \Bigg].
\end{equation}
We observe that, for every $\delta> 0$,
\begin{equation} \label{stimavalida}
\| \varphi(t) \|^2_H = \langle \varphi(t), \varphi(t) \rangle_{V^*, V} \leq \| \varphi(t) \|_{V^*}\| \varphi(t) \|_{V} 
\leq \frac{\delta}{2}\| \varphi(t) \|^2_{V} + \frac{1}{2\delta} \| \varphi(t) \|^2_{V^*}.
\end{equation} 
Choosing $\delta = \frac{b^2}{K + b^2}$ in \eqref{stimavalida}, we can estimate the first term of the right-hand side of \eqref{precedentissima3} as follows:
\begin{equation} \label{pezzoaggiunto}
(K + b^2)  \int_0^t \|\varphi(s)\|_H^2 \ ds  
\leq \frac{b^2}{2} \int_0^t \| \varphi(s) \|^2_{V} \ ds + \frac{(K+b^2)^2  \nu}{b^4 } \int_0^t \frac{b^2}{2 \nu}\| \varphi(s) \|^2_{V^*} \ ds.
\end{equation} 
Then, due to \eqref{pezzoaggiunto}, from \eqref{precedentissima3} we obtain that
\begin{equation} \nonumber
\frac{1}{2} \|\kappa(t)\|_H^2 + \int_{Q_t}  |\nabla \kappa|^2 +  \frac{b^2}{2} \int_0^t \| \varphi(s)\|_V^2 \ ds 
+ \frac{b^2}{2 \nu} \| \varphi(t) \|^2_{V^*}
\end{equation}
\begin{equation} \label{utilissima11112222245678}
\leq c \int_0^t \Bigg( \frac{1}{2}  \|\kappa(s)\|_H^2 +  \frac{b^2}{2 \nu}\| \varphi(s) \|^2_{V^*} \Bigg) \ ds
+ \frac{b^2}{2 \nu} \| \varphi_0 \|^2_{V^*}
+  \frac{1}{2} \|\kappa_0\|_H^2
+ 2 a^2 \|f\|_{L^2(0,T; H)}^2 + 3T \| \kappa^*\|_W^2.
\end{equation}
Due to  \eqref{f}--\eqref{regolaritaesistenza3},
the last four terms on the right-hand side of \eqref{utilissima11112222245678} are bounded uniformly with respect to $\varepsilon$ . Then,
by applying the Gronwall lemma, we conclude that
\begin{equation} \nonumber
 \|\kappa(t)\|_H + \|\nabla \kappa\|_{L^2(0,T;H)} +  \| \varphi\|_{L^2(0,T;V)}+  \| \varphi(t) \|_{V^*}
\end{equation}
\begin{equation} \label{precedentissima51}
\leq  c \bigg( \| \varphi_0 \|_{V^*}+  \|\kappa_0\|_H + \|f\|_{L^2(0,T; H)} +  \| \kappa^*\|_W \bigg)  
\end{equation} 
for some positive constant $c$ which depends only on $\Omega$, $T$
and the structure \eqref{beta}--\eqref{pi}, \eqref{A1}--\eqref{stimaA} and \eqref{parametri}--\eqref{0} of the system.
Now, we recall that \eqref{precedentissima51} is equivalent to
\begin{equation} \nonumber
\|\kappa_1 - \kappa_2\|_{L^{\infty}(0,T;H) \cap L^2(0,T;V)  } + \| \varphi_1 - \varphi_2\|_{L^{\infty}(0,T;V^*) \cap L^2(0,T;V) }
\end{equation} 
\begin{equation} \label{riferimentodipendenzacontinuafinale}
\leq  c \big( \| \varphi_{0_1} - \varphi_{0_2} \|_{V^*} + \|\kappa_{0_1} - \kappa_{0_2}\|_H + \|f_1 - f_2\|_{L^2(0,T; H)} +  \| \kappa^*_1- \kappa^*_2\|_W \big)   .
\end{equation} 
If $f_1 = f_2$, $\kappa^*_1 = \kappa^*_2$, $\kappa_{0_1} = \kappa_{0_2}$ and $\varphi_{0_1} = \varphi_{0_2}$, 
from \eqref{riferimentodipendenzacontinuafinale} we conclude that 
$\kappa_1 = \kappa_2$ and $\varphi_1 = \varphi_2$, i.e., the solution of problem $(\widetilde{P})$ (see \eqref{iniziale1-uni1}--\eqref{iniziale7-uni1}) is unique. From this fact, we immediately infer the uniqueness of the solution for our initial Problem $(P)$ (see \eqref{iniziale1}--\eqref{iniziale7}).

\setcounter{equation}{0}
\section{Sliding mode control}
This section is devoted to the proof of Theorem \ref{Teo3-Sliding-mode}. 
The argument we use in the proof relies in the following
Lemma (see \cite[Lemma 4.1, p. 20]{BaCoGiMaRo}).

\begin{Lemma-Sliding} \label{lemmautileslide}
Let $a_0$,  $b_0$, $\Upsilon_0$,  $\rho$ $\in \mathbb{R}$ be such that
\begin{equation} 
a_0, \ b_0,\ \Upsilon_0 \geq 0 \quad \textrm{and} \quad \rho > a_0^2 + 2 b_0 + 2 \frac{\Upsilon_0}{T}
\end{equation} 
and let $\Upsilon: \ [0,T] \rightarrow [0, + \infty)$ be an absolute ly continuous function satisfying $\Upsilon(0)=\Upsilon_0$
and
\begin{equation} 
\Upsilon' + \rho \leq a_0 \rho^{1/2} + b_0 \ \ \textrm{a.e. in the set $P:= \{ t\in (0,T): \Upsilon(t) > 0 \}$}.
\end{equation} 
The n, the  following conditions hold true :
\begin{enumerate}
	\item If $\Upsilon_0=0$, then $\Upsilon$ vanishe s ide ntically.
	\item If $\Upsilon_0>0$, the n the re  e xists $T^* \in (0,T)$ satisfying $T^* \leq 2 \Upsilon_0 /(\rho - a_0^2 - 2 b_0)$ such that
	$\Upsilon$ is strictly decreasing in $(0,T^* )$ and 	$\Upsilon$ vanishes in $[T^*,T]$.
\end{enumerate}
\end{Lemma-Sliding}

We assume $a=1$, $b= \ell$ and $A = \rho\Sign$ and 
consider the approximating problem ($\widetilde{P}_{\varepsilon}$)
obtained from ($P_{\varepsilon}$) (see \eqref{iniziale1-epsilon}--\eqref{iniziale7-epsilon}) 
with the usual change of variables 
\begin{equation} \label{camio222}
\kappa_{\varepsilon} = \chi_{\varepsilon} + \ell \varphi_{\varepsilon} - \kappa^* , 
\quad \quad \kappa_{0\varepsilon} = \chi_{0\varepsilon} + \ell \varphi_{0\varepsilon} - \kappa^*. 
\end{equation}  
We have that
\begin{equation} \label{iniziale1-uni1-rho}
\partial_t\kappa_{\varepsilon} - \Delta \kappa_{\varepsilon} 
+ \ell \Delta \varphi_{\varepsilon} - \Delta \kappa^* + \rho \sigma_{\varepsilon} = f_{\varepsilon} \ \ \textrm{a.e. in $Q$,} 
\end{equation}
\begin{equation} \label{iniziale2-uni1-rho}
\partial_t \varphi_{\varepsilon} - \Delta \mu_{\varepsilon}= 0  \ \ \textrm{a.e. in $Q$,}
\end{equation}
\begin{equation} \label{iniziale3-uni1-rho}
\mu_{\varepsilon}= - \nu \Delta \varphi_{\varepsilon} + \xi_{\varepsilon} + \pi(\varphi_{\varepsilon}) - \gamma (\kappa_{\varepsilon} - \ell \varphi_{\varepsilon} + \kappa^*) \ \ \textrm{a.e. in $Q$,} 
\end{equation}
\begin{equation} \label{iniziale4-uni1-rho}
\sigma_{\varepsilon}(t) \in \Signe(\kappa_{\varepsilon}(t)) \ \textrm{for a.e. $t \in (0,T)$,}  
\end{equation}
\begin{equation} \label{iniziale5-uni1-rho}
\xi_{\varepsilon} \in \beta_{\varepsilon}(\varphi_{\varepsilon}) \ \textrm{a.e. in $Q$,}
\end{equation}
\begin{equation} \label{iniziale6-uni1-rho}
\partial_{\bf{n}} \kappa_{\varepsilon} = \partial_{\bf{n}} \varphi_{\varepsilon} = \partial_{\bf{n}} \mu_{\varepsilon} = 0 \ \ \textrm{on $\Sigma$,}
\end{equation}
\begin{equation} \label{iniziale7-uni1-rho}
\kappa_{\varepsilon}(0) = \kappa_{0\varepsilon}, \quad \quad \varphi_{\varepsilon}(0) = \varphi_{0\varepsilon} \quad \textrm{in $\Omega$}.
\end{equation}

\paragraph{Further a priori uniform estimates.}
We test \eqref{iniziale1-uni1-rho} by $\partial_t\kappa_{\varepsilon}$ and integrate over $Q_t$.
Recalling that
\begin{equation}
\int_{Q_t} \rho \sigma_{\varepsilon} \partial_t\kappa_{\varepsilon}
= \rho \| \kappa_{\varepsilon} (t) \|_{H, \varepsilon} - \rho \| \kappa_0 \|_{H, \varepsilon},
\end{equation}
we have that
\begin{equation} \nonumber
\int_{Q_t} |\partial_t\kappa_{\varepsilon}|^2
+ \frac{1}{2} \int_{\Omega} |\nabla \kappa_{\varepsilon}(t)|^2 
+ \rho \| \kappa_{\varepsilon} (t) \|_{H, \varepsilon} 
= \frac{1}{2} \int_{\Omega} |\nabla \kappa_0|^2 
\end{equation}
\begin{equation} \label{slideccoultima222}
+ \rho \| \kappa_0 \|_{H, \varepsilon}
+ \int_{Q_t} \Delta \kappa^* \partial_t\kappa_{\varepsilon} 
+ \int_{Q_t} f_{\varepsilon} \partial_t\kappa_{\varepsilon}
- \int_{Q_t} \ell \Delta \varphi_{\varepsilon} \partial_t\kappa_{\varepsilon} .
\end{equation} 
We observe that $\| \kappa_0 \|_{H, \varepsilon} \leq \| \kappa_0 \|_H$ (cf. \eqref{segno1}).
Then, thanks to \eqref{0} and \eqref{tecnicadistima1}, the first two terms on the right-hand side of \eqref{slideccoultima222} 
are estimated as follows: 
\begin{equation} \label{aiuto1}
\frac{1}{2} \int_{\Omega} |\nabla \kappa_0|^2 + \rho \| \kappa_0 \|_{H, \varepsilon} \leq c (1 + \rho).
\end{equation} 
Due to \eqref{0} and \eqref{convergenzadif},
applying \eqref{dis2} to  the third and fo urth term o n the right-hand side o f \eqref{slideccoultima222},
we have that
\begin{equation} \label{aiuto2}
\int_{Q_t} \Delta \kappa^* \partial_t\kappa_{\varepsilon} 
\leq \frac{1}{4} \int_{Q_t} | \partial_t\kappa_{\varepsilon} |^2 + \int_{Q_t} |\Delta \kappa^*|^2 
= \frac{1}{4} \int_{Q_t} | \partial_t\kappa_{\varepsilon} |^2 + c,
\end{equation} 
\begin{equation} \label{aiuto3}
\int_{Q_t} f_{\varepsilon} \partial_t\kappa_{\varepsilon}
\leq \frac{1}{4} \int_{Q_t} |\partial_t\kappa_{\varepsilon}|^2 + \int_{Q_t} |f_{\varepsilon}|^2
\leq \frac{1}{4} \int_{Q_t} |\partial_t\kappa_{\varepsilon}|^2 + c.
\end{equation} 
Mo reo ver, integrating by parts the last term o f \eqref{slideccoultima222}, we formally have that
\begin{equation} \nonumber
- \int_{Q_t} \ell \Delta \varphi_{\varepsilon} \partial_t\kappa_{\varepsilon} 
= \ell \int_{Q_t} \nabla \varphi_{\varepsilon} \cdot \nabla(\partial_t\kappa_{\varepsilon} )
\end{equation}
\begin{equation} \label{paaarti2} 
= \ell \int_{\Omega} \nabla \varphi_{\varepsilon}(t) \cdot \nabla \kappa_{\varepsilon} (t)
- \ell \int_{\Omega} \nabla \varphi_0 \cdot \nabla \kappa_0
- \ell \int_{Q_t}  \nabla(\partial_t \varphi_{\varepsilon}) \cdot \nabla \kappa_{\varepsilon} .
\end{equation}
Using \eqref{dis2} and the H\"older inequality,
the first term o n the right-hand side o f \eqref{paaarti2} is estimated as fo llo ws:
\begin{equation} \nonumber
\Bigg| \ell \int_{\Omega} \nabla \varphi_{\varepsilon}(t) \cdot \nabla \kappa_{\varepsilon} (t) \Bigg|
\leq \frac{1}{4} \int_{\Omega} |\nabla \kappa_{\varepsilon} (t) |^2 + \ell^2 \int_{\Omega} |\nabla \varphi_{\varepsilon}(t) |^2
\end{equation}
\begin{equation} \nonumber
= \frac{1}{4} \int_{\Omega} |\nabla \kappa_{\varepsilon} (t) |^2 
+ \ell^2 \int_{\Omega} \bigg| \nabla \bigg( \varphi_0 + \int_0^t \partial_t \varphi_{\varepsilon}(s) \ ds \bigg) \bigg|^2
\end{equation}
\begin{equation} \nonumber
\leq
\frac{1}{4} \int_{\Omega} |\nabla \kappa_{\varepsilon} (t) |^2 
+ 2\ell^2 \int_{\Omega} |\nabla \varphi_0|^2 
+ 2\ell^2 \int_{\Omega} \bigg|  \int_0^t \nabla (\partial_t \varphi_{\varepsilon}(s)) \ ds \bigg|^2
\end{equation}
\begin{equation} \label{paaarti3}
\leq
\frac{1}{4} \int_{\Omega} |\nabla \kappa_{\varepsilon} (t) |^2 
+ 2\ell^2 \int_{\Omega} |\nabla \varphi_0|^2 
+ 2 T \ell^2 \int_{Q_t} | \nabla (\partial_t \varphi_{\varepsilon}) |^2 .
\end{equation}
Due to  \eqref{0}, the seco nd term o n the right-hand side o f \eqref{paaarti2} and similarly
the seco nd term o n the right-hand side o f \eqref{paaarti3}
are estimated by a po sitive co nstant $c$ independent o f $\rho$ and $\varepsilon$. Indeed
\begin{equation}
- \ell \int_{\Omega} \nabla \varphi_0 \cdot \nabla \kappa_0 \leq \ell^2 \int_{\Omega} |\nabla \varphi_0|^2 
+ \frac{1}{4} \int_{\Omega} |\nabla \kappa_0|^2 \leq c.
\end{equation}
Applying inequality \eqref{dis2} to the last term on the right-hand side of \eqref{paaarti2} we obtain that
\begin{equation} \label{aiuto4}
- \ell \int_{Q_t}  \nabla(\partial_t \varphi_{\varepsilon}) \cdot \nabla \kappa_{\varepsilon} 
\leq \frac{1}{4} \int_{Q_t} |\nabla \kappa_{\varepsilon} |^2 + \ell^2 \int_{Q_t} |\nabla(\partial_t \varphi_{\varepsilon})|^2.
\end{equation}
Then, thanks to \eqref{aiuto1}--\eqref{aiuto4}, from \eqref{slideccoultima222} we infer that
\begin{equation} \nonumber
\frac{1}{2} \int_{Q_t} |\partial_t\kappa_{\varepsilon}|^2
+ \frac{1}{4} \int_{\Omega} |\nabla \kappa_{\varepsilon}(t)|^2 
+ \rho \| \kappa_{\varepsilon} (t) \|_{H, \varepsilon} 
\end{equation} 
\begin{equation} \label{superaiuto1}
\leq c (1 + \rho) + \ell^2(1 + 2T) \int_{Q_t}  |\nabla (\partial_t \varphi_{\varepsilon})|^2
+ \frac{1}{4} \int_{Q_t}  |\nabla \kappa_{\varepsilon}|^2 .
\end{equation} 
Now, we formally differentiate \eqref{iniziale2-uni1-rho}
and \eqref{iniziale3-uni1-rho} with respect to time and obtain that
\begin{equation} \label{doppiainiziale2-uni1-rho}
\partial_{tt} \varphi_{\varepsilon} - \Delta \partial_{t}\mu_{\varepsilon}= 0,
\end{equation}
\begin{equation} \label{doppiainiziale3-uni1-rho}
\partial_t \mu_{\varepsilon}= - \nu \Delta \partial_t \varphi_{\varepsilon} 
+ \beta'_{\varepsilon}(\varphi_{\varepsilon}) \partial_t \varphi_{\varepsilon} 
+ \pi'(\varphi_{\varepsilon})\partial_t \varphi_{\varepsilon} 
- \gamma (\partial_t \kappa_{\varepsilon} - \ell \partial_t \varphi_{\varepsilon}).
\end{equation}
Acco rding to \eqref{eluiilveroriferimentoresatto}, $m(\partial_t \varphi_{\varepsilon})=0$. Co nsequently,
$\partial_t \varphi_{\varepsilon} \in D(\mathcal{N})$ and we can test \eqref{doppiainiziale2-uni1-rho} 
by $ \mathcal{N} (\partial_t \varphi_{\varepsilon})$
and \eqref{doppiainiziale3-uni1-rho} by $\partial_t \varphi_{\varepsilon}$, respectively.
Integrating over $Q_t$, we have that
\begin{equation} \label{nuovaRegg2}
- \int_{Q_t} \partial_t \mu_{\varepsilon} \partial_t \varphi_{\varepsilon} 
= \frac{1}{2}  \| \partial_t \varphi_{\varepsilon} (t) \|^2_{V^*}  
- \frac{1}{2}  \| \partial_t \varphi_{\varepsilon}(0) \|^2_{V^*}   ,
\end{equation}
\begin{equation} \nonumber
\int_{Q_t} \partial_t \mu_{\varepsilon} \partial_t \varphi_{\varepsilon}
= \nu \int_{Q_t} |\nabla \partial_t \varphi_{\varepsilon}|^2 
+\int_{Q_t} \beta_{\varepsilon}'(\varphi_{\varepsilon}) |\partial_t \varphi_{\varepsilon}|^2 
\end{equation} 
\begin{equation} \label{nuovaRegg3}
+ \int_{Q_t} \pi'(\varphi_{\varepsilon}) |\partial_t \varphi_{\varepsilon}|^2 
- \gamma \int_{Q_t} \partial_t \varphi_{\varepsilon} \partial_t\kappa_{\varepsilon}
+ \ell \gamma \int_{Q_t} |\partial_t \varphi_{\varepsilon}|^2.
\end{equation}
Combining \eqref{nuovaRegg2} and \eqref{nuovaRegg3} we obtain that
\begin{equation} \nonumber
\frac{1}{2} \| \partial_t \varphi_{\varepsilon} (t) \|^2_{V^*}
+ \nu \int_{Q_t} |\nabla \partial_t \varphi_{\varepsilon}|^2 
+ \ell \gamma \int_{Q_t} |\partial_t \varphi_{\varepsilon}|^2 
= \frac{1}{2}  \| \partial_t \varphi_{\varepsilon}(0) \|^2_{V^*}   
\end{equation}
\begin{equation} \label{slideccoultima}
- \int_{Q_t} \beta_{\varepsilon}'(\varphi_{\varepsilon}) |\partial_t \varphi_{\varepsilon}|^2 
- \int_{Q_t} \pi'(\varphi_{\varepsilon}) |\partial_t \varphi_{\varepsilon}|^2 
+ \gamma \int_{Q_t} \partial_t \kappa_{\varepsilon} \partial_t \varphi_{\varepsilon} .
\end{equation} 
Thanks to \eqref{0} and \eqref{tecnicadistima1},
the first term o n the right-hand side o f \eqref{slideccoultima} is bo unded by a po sitive co nstant
$c$ independent o f $\rho$ and $\varepsilon$ (cf. the analo go us bo und discussed belo w \eqref{utilissima111122222}). 
Since $\beta_{\varepsilon}$ is maximal mo no to ne, the seco nd term o n the right-hand side o f \eqref{slideccoultima} is no n-po sitive.
As $\pi$ is a Lipschitz co ntinuo us functio n with Lipschitz co nstant $C_{\pi}$, we have that
\begin{equation} \label{paaarti1}
- \int_{Q_t} \pi'(\varphi_{\varepsilon}) |\partial_t \varphi_{\varepsilon}|^2 
\leq C_{\pi} \int_{Q_t} |\partial_t \varphi_{\varepsilon}|^2.
\end{equation}
Finally, using \eqref{dis2}, the last term on the right-hand side of \eqref{slideccoultima} is
estimated as follows:
\begin{equation} \label{paaarti1333}
\gamma \int_{Q_t} \partial_t \kappa_{\varepsilon} \partial_t \varphi_{\varepsilon} 
\leq \frac{1}{4}\bigg( \frac{\nu}{\ell^2 (1+2T) +1} \bigg) \int_{Q_t} |\partial_t \kappa_{\varepsilon}|^2 
+  \gamma^2 \frac{\ell^2 (1+2T) +1}{\nu}  \int_{Q_t} |\partial_t \varphi_{\varepsilon}|^2,
\end{equation}
where the reaso n o f such invo lved co nstants will be clear in a mo ment.
Due to  \eqref{paaarti1}--\eqref{paaarti1333} and the previo us o bservatio ns, 
fro m \eqref{slideccoultima} we infer that
\begin{equation} \nonumber
\frac{1}{2} \| \partial_t \varphi_{\varepsilon} (t) \|^2_{V^*}
+ \nu \int_{Q_t} |\nabla \partial_t \varphi_{\varepsilon}|^2 
+ \ell \gamma \int_{Q_t} |\partial_t \varphi_{\varepsilon}|^2 
\end{equation}
\begin{equation} \label{slideccoultimaff2333}
\leq c + \frac{1}{4} \bigg( \frac{\nu}{\ell^2 (1+2T) +1} \bigg)\int_{Q_t} |\partial_t \kappa_{\varepsilon}|^2 
+  \bigg(\gamma^2 \frac{\ell^2 (1+2T) +1}{\nu}  + C_{\pi} \bigg) \int_{Q_t} |\partial_t \varphi_{\varepsilon}|^2.
\end{equation} 
Multiplying \eqref{slideccoultimaff2333} by $(\ell^2 (1+2T) +1)/ \nu$
and adding it to \eqref{superaiuto1}, we infer that
\begin{equation} \nonumber
\frac{1}{4} \int_{Q_t} |\partial_t\kappa_{\varepsilon}|^2
+ \frac{1}{4} \int_{\Omega} |\nabla \kappa_{\varepsilon}(t)|^2 
+ \rho \| \kappa_{\varepsilon} (t) \|_{H, \varepsilon} 
+ \int_{Q_t} |\nabla \partial_t \varphi_{\varepsilon}|^2 
+ C_1 \int_{Q_t} |\partial_t \varphi_{\varepsilon}|^2 
\end{equation} 
\begin{equation} \label{slideccoultimaff233344}
+ C_2 \| \partial_t \varphi_{\varepsilon} (t) \|^2_{V^*}
\leq c (1 + \rho) 
+ \frac{1}{4} \int_{Q_t}  |\nabla \kappa_{\varepsilon}|^2
+  C_3 \int_{Q_t} |\partial_t \varphi_{\varepsilon}|^2,
\end{equation} 
where
\begin{equation} \nonumber
C_1 = \frac{\ell^3 \gamma (1+2T) + \ell \gamma}{\nu}, \quad \quad
C_2 = \frac{\ell^2 (1+2T) +1}{2 \nu}, 
\end{equation} 
\begin{equation} \nonumber
C_3 = \gamma^2 \bigg(\frac{\ell^2 (1+2T) +1}{\nu}\bigg)^2  + C_{\pi} \frac{\ell^2 (1+2T) +1}{\nu} +\ell^2 (1+2T).
\end{equation} 
Denoting by $C_4$ the minimum between $1$ and $C_1$,
and applying the inequality \eqref{stimasimon1} with $\delta = \sqrt{C_4} / \sqrt{2C_3}$ 
to the last term on the right-hand side of  \eqref{slideccoultimaff233344}, we obtain that
\begin{equation}\label{helppp1}
C_3 \int_{Q_t} |\partial_t \varphi_{\varepsilon}|^2 
\leq \frac{C_4}{2} \int_0^t \|\partial_t \varphi_{\varepsilon} (s) \|_V^2 \ ds
+2 K^2 C_3 \int_0^t \| \partial_t \varphi_{\varepsilon} (s) \|_{V^*}^2 \ ds. 
\end{equation}
Thanks to \eqref{helppp1}, from \eqref{slideccoultimaff233344} we infer that
\begin{equation} \nonumber
\frac{1}{4} \int_{Q_t} |\partial_t\kappa_{\varepsilon}|^2
+ \frac{1}{4} \int_{\Omega} |\nabla \kappa_{\varepsilon}(t)|^2 
+ \rho \| \kappa_{\varepsilon} (t) \|_{H, \varepsilon} 
+ \frac{C_4}{2} \int_0^t \| \partial_t \varphi_{\varepsilon} (s) \|_V^2 
+ C_2 \| \partial_t \varphi_{\varepsilon} (t) \|^2_{V^*}
\end{equation} 
\begin{equation} \label{slideccoultimaff23334444}
\leq c (1 + \rho) 
+ \frac{1}{4} \int_{Q_t}  |\nabla \kappa_{\varepsilon}|^2
+2 K^2 C_3 \int_0^t \| \partial_t \varphi_{\varepsilon} (s) \|_{V^*}^2.
\end{equation} 
From \eqref{slideccoultimaff23334444}, by applying the Gronwall lemma, we conclude that
\begin{equation}
\| \partial_t\kappa_{\varepsilon}  \|_{L^2(0,T; H)}
+ \| \kappa_{\varepsilon}  \|_{L^{\infty}(0,T; V)}
+ \| \partial_t \varphi_{\varepsilon}  \|_{L^{\infty}(0,T; V^*)}
+ \| \partial_t \varphi_{\varepsilon}  \|_{L^2(0,T; V)} \leq c(1+ \rho^{1/2}),
\end{equation}
whence
\begin{equation} \label{intermezzo1}
\| \kappa_{\varepsilon} \|_{H^1(0,T;H) \cap L^\infty (0,T;V)} \leq c(1+ \rho^{1/2}),
\end{equation}
\begin{equation} \label{intermezzo2}
\| \varphi_{\varepsilon} \|_{W^{1, \infty}(0,T;V^*) \cap H^1(0,T;V)} \leq c(1+ \rho^{1/2}).
\end{equation}
Due to \eqref{intermezzo1}--\eqref{intermezzo2} and the change of variables stated by \eqref{camio222},
we have that
\begin{equation} \label{intermezzo3serio}
\| \chi_{\varepsilon} \|_{H^1(0,T;H) \cap L^\infty (0,T;V)} \leq c(1+ \rho^{1/2}).
\end{equation}
Proceeding as in the second a priori estimate 
(cf. \eqref{astrodiriferimento}--\eqref{stimagilardi})
and recalling \eqref{intermezzo2}--\eqref{intermezzo3serio}, 
from \eqref{ultimalemma2222} we infer that
\begin{equation} \label{intermezzo4serio}
\| \xi_{\varepsilon} \|_{L^\infty (0,T;L^1(\Omega))} \leq c(1+ \rho^{1/2}).
\end{equation}
No w, with the analo go us technique applied in the third a prio ri estimate,
thanks to  \eqref{intermezzo2}--\eqref{intermezzo4serio},
from \eqref{finalelemma3}  we obtain that
\begin{equation} \label{intermezzo5serio}
\| m(\mu_{\varepsilon}) \|_{L^\infty (0,T)} \leq c(1+ \rho^{1/2}).
\end{equation}
Then, due to \eqref{intermezzo5serio} and the Poincar$\acute{e}$ inequality,
by comparison in \eqref{iniziale2-uni1-rho} we deduce that
\begin{equation} \label{intermezzo6serio}
\| \mu_{\varepsilon} \|_{L^\infty (0,T;V)} \leq c(1+ \rho^{1/2}).
\end{equation}
Finally, with the same co mputatio ns as explained in the fifth a prio ri estimate 
(cf. \eqref{primalemma5}--\eqref{ultimosforzo}), thanks to  \eqref{intermezzo2}--\eqref{intermezzo3serio}
and \eqref{intermezzo6serio}, from \eqref{passaggiochiave5} we infer that
\begin{equation} \label{intermezzo7serio}
\| \xi_{\varepsilon} \|_{L^\infty (0,T;H)} \leq c(1+ \rho^{1/2}),
\end{equation}
whence, by comparison of every term in \eqref{iniziale3-uni1-rho}, we conclude that
\begin{equation} \label{intermezzo8serio}
\| \Delta \varphi_{\varepsilon} \|_{L^\infty (0,T;H)} \leq c(1+ \rho^{1/2}).
\end{equation}

\paragraph{Existence of sliding mode.}
Due to  \eqref{0}, \eqref{convergenzadif} and \eqref{intermezzo8serio},
we can rewrite \eqref{iniziale1-uni1-rho} in the form
\begin{equation} \label{sliding}
\partial_t\kappa_{\varepsilon}  - \Delta \kappa_{\varepsilon} + \rho\sigma_{\varepsilon}
= g_{\varepsilon} := f_{\varepsilon} - \ell \Delta \varphi_{\varepsilon} + \Delta \kappa^*,
\end{equation}
with
\begin{equation} \label{slidingstima}
\| g_{\varepsilon} \|_{L^{\infty}(0,T; H)} \leq c(1+ \rho^{1/2}),
\end{equation}
where $c$ depends o nly o n the structure and the data invo lved in the statement.
In o rder to  pro ve the existence o f sliding mo de, we fix the co nstant $c$  appearing in \eqref{slidingstima} and set
\begin{equation} \label{legamerot}
\rho^* := c^2 + 2c  + \frac{2}{T} \| \chi_0 + \ell \varphi_0 - \kappa^*\|_H
\end{equation}
and assume $\rho>\rho^*$. We also set
\begin{equation}
\Upsilon_{\varepsilon}(t) := \| \kappa_{\varepsilon}(t) \|_H \quad \textrm{for $t \in [0,T]$}.
\end{equation}
By assuming $h \in (0,T)$ and $t \in (0, T-h)$, we multiply \eqref{sliding} 
by $\sigma_{\varepsilon} = \Sign_{\varepsilon}(\kappa_{\varepsilon})$
and integrate o ver 
$(t,t+h) \times \Omega$. We have that
\begin{equation} \nonumber
\int_t^{t+h} (\partial_t\kappa_{\varepsilon}(s), \sigma_{\varepsilon}(s))_H \ ds 
+  \int_t^{t+h} \int_{\Omega} \nabla \kappa_{\varepsilon} \cdot \nabla \sigma_{\varepsilon}
+ \rho \int_t^{t+h} \| \sigma_{\varepsilon} (s) \|^2_H \ ds
\end{equation}
\begin{equation} \label{slide2}
= \int_t^{t+h} (g_{\varepsilon}(s) , \sigma_{\varepsilon}(s))_H \ ds .
\end{equation}
Recalling that $\Sign_{\varepsilon}(v)$ is the gradient at $v$ o f the $C^1$ functio nal $\| \cdot \|_{H, \varepsilon}$,
fro m \eqref{segno1}--\eqref{segno2} we deduce that
\begin{equation} \nonumber
(\partial_t\kappa_{\varepsilon}(s), \sigma_{\varepsilon}(s))_H = \frac{d}{dt} \int_0^{\Upsilon_{\varepsilon}(t)} \min{\{ s/\varepsilon, \ 1 \}} \ ds \quad
\textrm{for a.a. $t \in (0,T)$.}
\end{equation}
Then, fo r the first term o n the right-hand side o f \eqref{slide2} we have that
\begin{equation} \nonumber
\int_{t}^{t+h} (\partial_t\kappa_{\varepsilon}(s), \sigma_{\varepsilon}(s))_H \ ds 
= \int_{\Upsilon_{\varepsilon}(t)}^{\Upsilon_{\varepsilon}(t+h)} \min{ \{ s/\varepsilon,  1 \}} \ ds.
\end{equation}
We also  no tice that \eqref{segno2} implies that
\begin{equation} \nonumber
\nabla \kappa_{\varepsilon}(t) \cdot \nabla \sigma_{\varepsilon} (t) 
= \frac{|\nabla \kappa_{\varepsilon}(t)|^2}{ \max{ \{ \varepsilon,  \| \kappa_{\varepsilon}(t) \|_H \} }} \geq 0 \ \ \ \textrm{a.e. in $\Omega$, for a.e. $t \in (0,T)$},
\end{equation}
whence the second integral on the left-hand side o f \eqref{slide2}
is no nnegative. Mo reo ver, as $ \| \sigma_{\varepsilon} (s) \|_H \leq 1$ for every $s $ (see \eqref{Sign3333333})
and \eqref{slidingstima} holds, we infer from \eqref{slide2} that
\begin{equation} \label{soccorso2}
\int_{\Upsilon_{\varepsilon}(t)}^{\Upsilon_{\varepsilon}(t+h)} \min{ \{ s / \varepsilon,  1 \}} \ ds + \rho \int_t^{t+h} \| \sigma_{\varepsilon} (s) \|^2_H \ ds \leq hc(\rho^{1/2} +1).
\end{equation}
At this point, we let $\varepsilon \searrow 0 $.
Due to \eqref{limiteforte1}--\eqref{limiteforte2}, \eqref{camio222}
and the uniqueness of the solution of the limit Problem \eqref{iniziale1-uni1}--\eqref{iniziale7-uni1}
(cf. Theorem \ref{Teorema-unicita}) we have that
\begin{equation} 
\kappa_{\varepsilon}  \rightarrow \kappa  \quad \quad \textrm{in $C^0(0,T; H)$}.
\end{equation} 
Besides, using standard weak, weakstar and compactness results, from \eqref{soccorso2} we infer that
\begin{equation}
\sigma_{\varepsilon}    \rightharpoonup^*  \sigma      \quad \quad \textrm{in $ L^{\infty} (0,T;H)$}.
\end{equation}
Then, taking the limit as $\varepsilon \searrow 0$
in \eqref{soccorso2} and denoting by 
\begin{equation}
\Upsilon(t) := \| \kappa(t) \|_H \quad \quad \textrm{for $t \in [0,T]$},
\end{equation}
we obtain that
\begin{equation} \nonumber
\Upsilon(t+h)- \Upsilon(t) +  \rho \int_t^{t+h} \| \sigma (s) \|^2_H \ ds 
\end{equation}
\begin{equation} \label{ultimaslide}
\leq \lim_{\varepsilon \searrow 0} \int_{\Upsilon_{\varepsilon}(t)}^{\Upsilon_{\varepsilon}(t+h)} \min{ \{ s / \varepsilon,  1 \}} \ ds 
+ \rho \liminf_{\varepsilon \searrow 0} \int_t^{t+h} \| \sigma_{\varepsilon} (s) \|^2_H \ ds \leq hc(\rho^{1/2} +1)
\end{equation}
for every $h \in (0,T)$ and $t \in (0,T-h)$.
Finally, we multiply \eqref{ultimaslide} by $1 / h$ and let $h$ tend to zero. We conclude that
\begin{equation} 
\Upsilon'(t) + \rho \| \sigma (t) \|^2_H \leq c(\rho^{1/2} +1) \ \ \textrm{for a.a. $t \in (0,T)$}.
\end{equation}
As $\| \sigma (t) \|_H = 1$ if $\| \kappa (t) \|_H > 0$ (see \eqref{Sign3333333}), 
we can apply Lemma \ref{lemmautileslide} with $a_0=b_0=c$ and
we observe that our condition $\rho> \rho^*$ completely fits the assumptions by \eqref{legamerot}.
Thus, we find $T^* \in [0,T)$ such that $\kappa(t) = 0 $ for every $t \in [T^*, T]$, i.e., \eqref{finedellarticolo3}.

\section*{Acknowledgments}
The author is very grateful to Professor Pierluigi Colli
for his advice, for his kind helpfulness and for several precious discussions.
Moreover, MC expresses his gratitude to
the referee for the careful reading of the manuscript 
and for a number of useful suggestions. Some
partial support from the GNAMPA (Gruppo Nazionale per l'Analisi Matematica, 
la Pro\-babilit\`a e le loro Applicazioni) of INdAM is acknowledged.


\begin{thebibliography}{99}

\bibitem{Bar10}
V.\ Barbu,
``Nonlinear differential equations of monotone types in Banach spaces'',
Springer, New York, 2010.
	
\bibitem{BaCoGiMaRo}
V.\ Barbu, P.\ Colli, G.\ Gilardi, G.\ Marinoschi and E.\ Rocca,
Sliding mode control
for a nonlinear phase-field system, SIAM J. Control Optim, to appear 2017
(see also preprint arXiv:1506.01665~[math.AP] (2015), 1--28).

\bibitem{Bre73}
H.\ Brezis,
``Op\'erateurs maximaux monotones et semi-groupes de contractions
dans les espaces de Hilbert'',
North-Holland Math. Stud. {\bf 5}, North-Holland, Amsterdam, 1973.

\bibitem{BrokSpr} 
M. \ Brokate and J.\ Sprekels,
``Hysteresis and phase transitions'',
Springer, New York, 1996.

\bibitem{Cag}
G.\ Caginalp,
An analysis of a phase field model of a free boundary,
Arch. Rational Mech. Anal. {\bf 92} (1986), 205--245.

\bibitem{CH58}
	J.\ W.\ {C}ahn and J.\ E.\ {H}illiard, 
	\newblock Free energy of a nonuniform system I. Interfacial free energy, 
	\newblock J.\ Chem.\ Phys. {\bf 2} (1958), 258--267.

\bibitem{eq1}
 L.\ Cherfils, A.\ Miranville and S.\ Zelik, The Cahn--Hilliard equation with logarithmic
potentials, Milan J. Math. \textbf{79} (2011), 561--596.

\bibitem{SC1}
P. \ Colli, M. \ H. \ Farshbaf-Shaker, G. \ Gilardi and J. \ Sprekels,
Optimal boundary control of a viscous Cahn–Hilliard system
with dynamic boundary condition and double obstacle potentials,
SIAM J. Control Optim.  \textbf{53} (2015), 2696--2721.

\bibitem{SC2}
P. \ Colli, M. \ H. \ Farshbaf-Shaker, G. \ Gilardi and J. \ Sprekels,
Second-order analysis of a boundary control problem for the viscous
Cahn–Hilliard equation with dynamic boundary condition,
Ann. Acad. Rom. Sci. Ser. Math. Appl. \textbf{7} (2015), 41--66.

\bibitem{CF15a} 
	P.\ {C}olli and T.\ {F}ukao, 
	\newblock {C}ahn--{H}illiard equation with dynamic boundary conditions and mass constraint on the boundary, 
	\newblock J.\ Math.\ Anal.\ Appl. {\bf 429} (2015), 1190--1213.

\bibitem{CF15bbbbbb} 
	P.\ {C}olli and T.\ {F}ukao, 
	\newblock Equation and dynamic boundary condition of {C}ahn--{H}illiard type with singular potentials, 
Nonlinear Anal. {\bf 127} (2015), 413--433.

\bibitem{CF15200009} 
P.\ {C}olli and T.\ {F}ukao, 
Nonlinear diffusion equations as asymptotic limits of Cahn–Hilliard systems,
J. Differential Equations {\bf 260} (2016), 6930--6959.

\bibitem{CGM}
P.\ Colli, G.\ Gilardi and G.\ Marinoschi,
A~boundary control problem for a possibly singular phase field system with dynamic boundary conditions,
J. Math. Anal. Appl. \textbf{434} (2016), 432--463.

\bibitem{CoGiMaRo}
P.\ Colli, G.\ Gilardi, G.\ Marinoschi and E.\ Rocca,
Optimal control for a phase field system with a possibly singular potential,
Math. Control Relat. Fields \textbf{6} (2016), 95--112.

\bibitem{eq2}
P.\ Colli, G.\ Gilardi and J.\ Sprekels, On the Cahn--Hilliard equation with dynamic
boundary conditions and a dominating boundary potential, J. Math. Anal. Appl. \textbf{419}
(2014), 972--994.

\bibitem{SC3}
P. \ Colli, G. \ Gilardi and J. \ Sprekels,
A boundary control problem for the pure Cahn–Hilliard equation with dynamic boundary conditions,
Adv. Nonlinear Anal. \textbf{4} (2015), 311--325.

\bibitem{SC4}
P. \ Colli, G. \ Gilardi and J. \ Sprekels,
A boundary control problem for the viscous Cahn–Hilliard equation with dynamic boundary conditions,
Appl. Math. Optim.  \textbf{73} (2016), 195--225.

\bibitem{collisprek}
P.\ Colli, G.\ Gilardi and J.\ Sprekels,
Constrained evolution for a quasilinear parabolic equation,
J. Optim. Theory Appl. \textbf{170} (2016), 713--734.

\bibitem{CLS99}
	P.\ {C}olli, Ph.\ {L}auren{\c{c}}ot and J.\ {S}prekels,
	\newblock Global solution to the {P}enrose--{F}ife phase field model with special heat flux laws,
	\newblock 181--188, in
	 ``Variations of domain and free-boundary problems in solid  mechanics'', 
	\newblock Solid Mech.\ Appl. {\bf 66},
	Kluwer Acad.\ Publ., Dordrecht, 1999.
	
\bibitem{CoGiMaRorrrrrrrrr}
P.\ Colli, G.\ Marinoschi and E.\ Rocca,
Sharp interface control in a Penrose-Fife model,
ESAIM Control Optim. Calc. Var. \textbf{22} (2016), 473--499.


\bibitem{Michele} 
	M.\ Colturato, 
	Solvability of a class of phase field systems related to a sliding mode control problem, 
	Appl. Math. \textbf{6} (2016), 623--650. 
	
\bibitem{EZ86}
	C.\ M.\ {E}lliott and S.\ {Z}heng, 
	\newblock On the {C}ahn--{H}illiard equation, 
	\newblock Arch.\ Rational\ Mech.\ Anal. {\bf 96} (1986), 339--357.
	
\bibitem{EllZheng}
{C.\ M.\ Elliott and S. Zheng,} 
Global existence and stability of solutions to the phase-field equations, 
Internat. Ser. Numer. Math. {\bf 95}, 46--58, in ``Free boundary problems'', Birkh\"auser
Verlag, Basel, 1990.
	

\bibitem{GMS09} 
	 G.\ {G}ilardi, A.\ {M}iranville and G.\ {S}chimperna, 
	\newblock On the {C}ahn--{H}illiard equation with irregular potentials and dynamic boundary conditions,
	\newblock Commun.\ Pure.\ Appl.\ Anal. {\bf 8} (2009), 881--912.

\bibitem{GraPetSch}
M.\ Grasselli, H.\ Petzeltov\'a and G.\ Schimperna,
Long time behavior of solutions to the Caginalp system with singular potential, 
Z. Anal. Anwend. {\bf 25} (2006), {51--72}.

\bibitem{preamboloiniziale}
M.\ Heida, Existence of solutions for two types of generalized versions of the Cahn--Hilliard equation, Appl. Math. \textbf{60} (2015), 51--90.

\bibitem{control1}	
M.\ Hinterm\"uller and D. Wegner, Distributed optimal control of the Cahn--Hilliard
system including the case of a double-obstacle homogeneous free energy density,
SIAM J. Control Optim. \textbf{50} (2012), 388--418.

\bibitem{control7}	
M.\ Hinterm\"uller and D.\ Wegner, Optimal control of a semi-discrete Cahn--Hilliard–
Navier–Stokes system, SIAM J. Control Optim. \textbf{52 } (2014), 747--772.
	
\bibitem{HKKY}
K.\ H.\ Hoffmann, N.\ Kenmochi, M.\ Kubo and N.\ Yamazaki, 
Optimal control problems for models of phase-field type with hysteresis of play operator,
 Adv. Math. Sci. Appl. {\bf 17} (2007), {305--336}.

\bibitem{I76}
U. \ Itkis, ``Control systems of variable structure'',
Wiley, Hoboken, 1976.

\bibitem{KenmNiez}
N.\ Kenmochi and M.\ Niezg\'odka,
Evolution systems of nonlinear variational inequalities arising from phase change problems, 
 Nonlinear Anal. {\bf 22} (1994), 1163--1180.

\bibitem{KN96}
	 N.\ Kenmochi and M.\ Niezg\'odka, 
	\newblock Viscosity approach to modelling non--isothermal diffusive phase separation,
	\newblock Japan.\ J.\ Indust.\ Appl.\ Math. {\bf 13} (1996), 135--169.

\bibitem{Kub12} 
	M.\ {K}ubo,
	\newblock The {C}ahn--{H}illiard equation with time-dependent constraint, 
	\newblock Nonlinear Anal. {\bf 75} (2012), 5672--5685. 

\bibitem{Lau}
Ph.\ Lauren\c cot,
Long-time behaviour for a model of phase-field type, 
 Proc. Roy. Soc. Edinburgh Sect.~A {\bf 126} (1996), 167--185.
	
\bibitem{O00}
Y.\ V.\ Orlov, 
Discontinuous unit feedback control of uncertain infinitedimensional
systems, 
 IEEE Trans. Automatic Control {\bf 45} (2000), 834--843.

\bibitem{eq3}
J.\ Pruss, R.\ Racke and S.\ Zheng, Maximal regularity and asymptotic behavior of
solutions for the Cahn--Hilliard equation with dynamic boundary conditions, Ann.
Mat. Pura Appl. (4) \textbf{185} (2006), 627--648.

\bibitem{eq4}
R.\ Racke and S.\ Zheng, The Cahn--Hilliard equation with dynamic boundary conditions,
Adv. Differential Equations \textbf{8} (2003), 83--110.

\bibitem{Show}
R. \ E. \ Showalter,
``Monotone operators in Banach Space and Nonlinear Partial Differential Equations'', 
AMS, Mathematica Surveys and Monographs {\bf 49},
1991.

\bibitem{Sim87}
	J.\ {S}imon, 
	\newblock Compact sets in the spaces $L^p(0,T;B)$, 
	\newblock Ann.\ Mat.\ Pura.\ Appl.~(4) {\bf 146} (1987), 65--96.


\bibitem{sprekzen}
J.\ Sprekels and S.\ M.\ Zheng,
Global smooth solutions to a thermodynamically consistent model of phase-field type in higher space dimensions,
J. Math. Anal. Appl. \textbf{176} (1993), 200--223.



\bibitem{control8}
	Q.\ F.\ Wang, Optimal distributed control of nonlinear Cahn--Hilliard systems with
computational realization, J. Math. Sci. (N. Y.) \textbf{177} (2011), 440--458.
	

	\bibitem{eq5}
H.\ Wu and S.\ Zheng, Convergence to equilibrium for the Cahn--Hilliard equation
with dynamic boundary conditions, J. Differential Equations \textbf{204} (2004), 511--531.
	
	\bibitem{control4}
X.\ P.\ Zhao and C.\ C.\ Liu, Optimal control of the convective Cahn--Hilliard equation,
Appl. Anal. \textbf{92} (2013), 1028--1045.

	\bibitem{control5}
X.\ P.\ Zhao and C.\ C.\ Liu, Optimal control for the convective Cahn--Hilliard equation
in 2D case, Appl. Math. Optim. \textbf{70} (2014), 61--82.


\end{thebibliography}
\end{document}